\documentclass[12pt]{amsart}

\input{xy}
\xyoption{all}
\usepackage{graphicx}
\usepackage{color}
\usepackage{verbatim}

\newtheorem{theorem}{Theorem} [section]
\newtheorem{prop}[theorem]{Proposition}
\newtheorem{lemma}[theorem]{Lemma}

\newtheorem{question}[theorem]{Question}

\theoremstyle{definition}

\newtheorem{example}[theorem]{Example}
\newtheorem{remark}[theorem]{Remark}

\numberwithin{equation}{section}
\numberwithin{figure}{section}

\newcommand\C{{\mathbb C}}

\renewcommand\P{{\mathbb P}}
\newcommand\bP{{\bf P}}
\newcommand\R{{\mathbb R}}
\newcommand\Z{{\mathbb Z}}
\newcommand\Q{{\mathbb Q}}

\renewcommand\phi{\varphi}
\renewcommand\O{\mathcal{O}} 
\newcommand{\lra}{\longrightarrow}

\newcommand\PSL{\mathrm{PSL}}

\newcommand\kbar{\bar{k}}

\DeclareMathOperator{\lhat}{{\widehat{\lambda}}}
\DeclareMathOperator{\hhat}{{\widehat{h}}}

\newcommand{\OO}{{\mathcal O}}
\newcommand{\N}{{\mathbb N}}
\newcommand{\bG}{{\mathbb G}}

\newcommand\iso{\simeq} 

\newcommand\Res {\operatorname{Res}} 

\newcommand\ord{\operatorname{ord}}  

\newcommand\M {\mathrm{M}} 
  
\newcommand\Mbar {\overline{\M}}

\begin{document}

\title{Rationality of dynamical canonical height}

\author{Laura De~Marco}
\address{
Laura DeMarco\\
Department of Mathematics\\
Northwestern University\\
2033 Sheridan Road\\
Evanston, IL 60208-2730 \\
USA
}
\email{demarco@northwestern.edu}

\author{Dragos Ghioca}
\address{
Dragos Ghioca\\
Department of Mathematics\\
University of British Columbia\\
Vancouver, BC V6T 1Z2\\
Canada
}
\email{dghioca@math.ubc.ca}

\date{\today}

\begin{abstract}
We present a dynamical proof of the well-known fact that the N\'eron-Tate canonical height (and its local counterpart) takes rational values at points of an elliptic curve over a function field $k = \C(X)$, where $X$ is a curve.  More generally, we investigate the mechanism for which the local canonical height for a map $f: \P^1\to \P^1$ defined over a function field $k$ can take irrational values (at points in a local completion of $k$), providing examples in all degrees $\deg f \geq 2$.  Building on Kiwi's classification of non-archimedean Julia sets for quadratic maps \cite{Kiwi:quad}, we give a complete answer in degree 2 characterizing the existence of points with irrational local canonical heights.  As an application we prove that if the heights $\hhat_f(a), \hhat_g(b)$ are rational and positive, for maps $f$ and $g$ of multiplicatively independent degrees and points $a, b \in \P^1(\kbar)$, then the orbits $\{f^n(a)\}_{n\geq 0}$ and $\{g^m(b)\}_{m\geq 0}$ intersect in at most finitely many points, complementing the results of \cite{Ghioca:Tucker:Zieve}.  
\end{abstract}



\maketitle

\thispagestyle{empty}

\section{Introduction}
\label{section introduction}

Throughout this article, we let $K$ be an algebraically closed field of characteristic $0$, $X$ a smooth projective curve defined over $K$, and $k=K(X)$ its function field.  We examine the values of the canonical height function associated to dynamical systems
	$$f: \P^1 \to \P^1$$
defined over $k$, of degree $\deg(f) \geq 2$.  The canonical height $\hhat_f$ was introduced in \cite{Call:Silverman} for general polarized dynamical systems $f$ on projective varieties; it measures arithmetic complexity in orbits.  In the special case where $f$ is the quotient of an endomorphism of an elliptic curve, the canonical height coincides with the N\'eron-Tate canonical height on the elliptic curve.  
For each place $v$ of the function field $k/K$, Call and Silverman showed that, in the presence of a ``weak N\'eron model," the values of the local canonical heights $\lhat_{f,v}$ are equal to certain intersection numbers (see \cite[Section~6]{Call:Silverman}), exactly as for abelian varieties.  Unfortunately, these weak N\'eron models often fail to exist in the dynamical setting (see \cite[Theorem~3.1]{Hsia:weak}).  Here, we take a more general approach to computing the local heights, building on the methods of \cite[Section 5]{Call:Silverman}, by analyzing the action of $f$ on the Berkovich projective line $\bP^1_v$ over a local field $\C_v$.  Our main results show that the local canonical heights can take irrational values, in contrast to the setting of elliptic curves and abelian varieties where all local and global height values are rational.

\subsection{Elliptic rationality}
Suppose that $E$ is an elliptic curve over $k = K(X)$; the \emph{N\'eron-Tate canonical height}  $\hhat_E:  E(\kbar) \to \R$ is defined as follows (see, e.g., \cite{Silverman:Advanced}). We let $\pi :E\lra \P^1$ be a nonconstant even function (such as the identification of a point $P$ with $-P$) and define
\begin{equation}
\label{definition canonical height elliptic curve} 
\hhat_E(P):=\frac{1}{\deg(\pi)}\cdot \lim_{m\to \infty} \frac{h(\pi([m]P))}{m^2},
\end{equation}
where $[m]$ is the multiplication-by-$m$ map on the elliptic curve $E$ and $h$ is a logarithmic Weil height on $\P^1(\kbar)$.  As shown in \cite[Chapter~III]{Silverman:Advanced}, $\hhat_E(P)\ge 0$ for each $P\in E(\kbar)$, with equality if and only if $P$ is a torsion point for $E$ or the pair $(E,P)$ is \emph{isotrivial} in the sense that $E$ is isomorphic (through some isomorphism $\iota$ defined over a finite extension of $k$) to an elliptic curve $E_0$ defined over $K$, while  $\iota(P)\in E_0(K)$.  The canonical height has a decomposition into N\'eron local heights,
	$$\hhat_E = \sum_{v\in M_k}  \lhat_{E,v},$$
where each $\lhat_{E,v}$ is well defined on the nonzero points of $E$ and $M_k$ denotes the set of all places of the function field $k/K$.   

For each place $v\in M_k$, we let $\C_v$ be a minimal complete and algebraically closed extension of $k$. 
We begin this article with a new, dynamical proof of a well-known ``local rationality" statement:

\begin{theorem}  \label{dynamical E}  
Suppose $E$ is an elliptic curve over $k$.  Fix any place $v \in M_k$.  Then, for any point $P \in E(\C_v)\setminus\{O\}$, the local canonical height $\lhat_{E,v}(P)$ is a rational number.   Consequently, for all $P \in E(\kbar)$ the canonical height $\hhat_E(P)$ is a rational number. 
\end{theorem}

\noindent
See \cite[Chapter~III,~Theorem~9.3]{Silverman:Advanced} for a proof that $\hhat_E(P)\in \Q$ for any $P\in E(\kbar)$; see \cite[Section 6 and p.~203]{Call:Silverman} and also \cite[Chapter~11~Theorem~5.1]{Lang:Diophantine} for the proofs that each local function $\lhat_{E,v}$ also takes values in $\Q$.  Both of these facts can be explained by means of intersection theory since  the local heights correspond to intersection numbers on a N\'eron local model.  

\subsection{Local irrationality}
As for elliptic curves, the canonical height of a morphism $f: \P^1\to \P^1$ of degree $d\geq 2$ has a local decomposition, 
	$$\hhat_f(x)=\sum_{v\in M_k} \lhat_{f,v}(x)$$
(see \S \ref{transforming}).  Also, as for elliptic curves, it is easy to show that local rationality holds for polynomial maps, as observed in \cite{Ingram:polyvariation} (see also the discussion immediately following \cite[Fact~3.6]{GT06}); that is, for any polynomial map $f: \P^1\to \P^1$ defined over $k$ and every $a\in \P^1(\C_v)\setminus\{\infty\}$, the value of $\lhat_{f,v}(a)$ is a rational number because of the strongly-attracting behavior of $\infty$.  Consequently, the global canonical height $\hhat_f(P)$ is rational for all $P \in \P^1(\kbar)$ and any polynomial $f$.  

A key ingredient in our proof of Theorem \ref{dynamical E} -- in studying the dynamics of an associated map on $\P^1$ -- is that all Type II points in the Julia set within the Berkovich projective line $\bP^1_v$ are preperiodic.  This is true for all rational maps $f$ defined over $k$ \cite[Proposition 3.9]{DF:degenerations2}, but nevertheless, the local heights are not always rational.  

%

\begin{theorem} \label{irrational}
For every degree $d\geq 2$, there exist a rational function $f \in k(z)$ of degree $d$, a place $v\in M_k$, and a point $x \in \P^1(\C_v)$ for which the local canonical height $\lhat_{f,v}(x)$ lies in $\R\setminus\Q$.   
\end{theorem}

The examples are inspired by the classical dynamical construction of a Markov partition, encoding orbits by their itineraries; see Proposition \ref{irrational general}.  In fact, there is nothing special about the irrational values that arise as local heights, as exhibited by the following example:

\begin{theorem} \label{all reals}
Let $k = \C(t)$ and $v = \ord_{t=0}$ so that $v(t) = 1$.  Let 
	$$f(z) = \frac{(z+1)(z-t)}{z+t}.$$
Then for all real numbers $-1 \leq \alpha \leq 0$, there exists a point $x \in \P^1(\C_v)$ so that 
	$$\lhat_{f,v}(x) = \alpha.$$
\end{theorem}

\begin{remark}
By choosing other normalizations for $\lhat_{f,v}$ in Theorem \ref{all reals}, the interval $[-1,0]$ of values can be shifted to any interval of the form $[n, n+1]$ with $n\in \Z$.  See \S\ref{properties}.
\end{remark}

Recall that polynomials only have rational values for their local heights.  For degree 2 maps, building on the work of Kiwi \cite{Kiwi:quad}, we give a complete description of which maps can produce irrational local heights.

\begin{theorem}  \label{quadratic}
Fix a place $v \in M_k$.  For each $f\in k(z)$ of degree $d=2$, exactly one of the following holds:
\begin{enumerate}
\item	$f$ has potential good reduction at $v$; 
\item $f$ is strongly polynomial-like at $v$; or
\item there are points in $\P^1(\C_v)$ with irrational local canonical heights.
\end{enumerate}
\end{theorem}

\noindent
The map $f$ is said to have \emph{potential good reduction} if, after a change of coordinates by a M\"obius transformation $\mu$ defined over $\kbar$, the reduction $\mu \circ f \circ \mu^{-1} \mod v$ has the same degree as $f$. For a map $f$ of degree $>1$, this will hold at all but finitely many places $v$ of $k$, and it holds at all places $v$ if and only if $f$ is isotrivial \cite{Baker:functionfields} (see also \cite{PST} when $f$ is an endomorphism of $\mathbb{P}^N$).  We say a quadratic map $f$ is \emph{strongly polynomial-like} at $v$ if its three fixed-point multipliers $\{\alpha, \beta, \gamma\}$ satisfy $|\alpha|_v = |\beta|_v > 1$ and $|\gamma|_v < 1$.  Note that the quadratic polynomials are characterized, up to conjugacy by an automorphism of $\P^1$, by the condition that one fixed-point multiplier is $\gamma =0$; furthermore, if the polynomial does not have potential good reduction, then the fixed point formula guarantees that the other two multipliers satisfy $|\alpha|_v = |\beta|_v > 1$.  (See \S\ref{multiplier}.)  

From \cite[Theorems 2 and 3]{Kiwi:quad}, our strong polynomial-like condition implies that the Julia set of $f$ is a Cantor set contained in $\P^1(\C_v)$, with $f|J(f)$ conjugate to the full one-sided 2-shift, but the conditions are not equivalent.  There exist maps with Cantor Julia set that do admit points with irrational local height; these are precisely the quadratic maps with fixed point multipliers satisfying $|\alpha|_v > |\beta|_v > 1 > |\gamma|_v = |1/\beta|_v > 0$.  Our proof of Theorem \ref{quadratic} shows that, in this polynomial-like but not strongly polynomial-like case, the points with irrational local height are necessarily in the Cantor Julia set.

\subsection{Global rationality?}
The constructions we exploit in Theorem \ref{irrational} suggest that the points $x\in \P^1(\C_v)$ with irrational local height will be transcendental over $k$.  We expect an affirmative answer to the following question.

\begin{question}
Given any rational function $f$ of degree $\geq 2$ defined over $k$, is $\hhat_f(P)$ a rational number for every $P \in \P^1(\kbar)$?
\end{question}

\subsection{Orbit intersections}  We conclude with an application of rationality, partly motivating our study. For any rational function $f$ of degree larger than $1$ defined over a field $F$ of characteristic $0$, and for any $a\in \P^1(F)$, we denote by $\OO_f(a)$ the orbit of $a$ under the iterates of $f$, i.e., the set of all $f^n(a)$ for $n\ge 0$.  It is natural to ask under what hypotheses the intersection $\OO_f(a)\cap \OO_g(b)$ can be infinite.  For example, if $f$ and $g$ share an iterate, and if $a$ and $b$ lie in the same grand orbit for this iterate, then the intersection is infinite.  It is expected that infinite intersection can only arise from such special arrangements. It is known that if $f$ and $g$ are polynomials, then an $f$-orbit can intersect a $g$-orbit in infinitely many points only when $f$ and $g$ share a common iterate (see \cite{GTZ} and \cite{Ghioca:Tucker:Zieve}).  On the other hand, almost nothing is known for general rational functions.   

The question was originally motivated by the Dynamical Mordell-Lang Conjecture (see \cite{GT-JNT}).  A special case of the Dynamical Mordell-Lang Conjecture predicts that if $C$ is an irreducible curve inside an arbitrary variety $X$ defined over a field $K$ of characteristic $0$ and endowed with an endomorphism $\Phi:X\longrightarrow X$, if there exists $\alpha\in X(K)$ such that $C$ contains infinitely many points in common with the orbit of $\alpha$ under $\Phi$, then $C$ must be periodic.  

We work over function fields $k/K$ as above, and we prove the result for rational functions under additional simplifying hypotheses:

\begin{theorem}
\label{main result application}
Fix $a, b\in \P^1(k)$ and $f,g\in k(z)$ with $d = \deg(f)$ and $e = \deg(g)$ both $>1$. Assume the following hypotheses are met: 
\begin{itemize}
\item[(1)] $d$ and $e$ are multiplicatively independent; 
\item[(2)]  the pair $(f,a)$ is not isotrivial over $k$;
\item[(3)]  $\hhat_f(a),\hhat_g(b)\in \Q$.
\end{itemize}
Then $\OO_f(a)\cap\OO_g(b)$ is finite.

Furthermore, if hypotheses~(2)-(3) hold, while condition~(1) is replaced by the weaker hypothesis:
\begin{itemize}
\item[(1')] $f$ and $g$ do not share a common iterate, 
\end{itemize}
then the Dynamical Mordell-Lang Conjecture yields that $\OO_f(a)\cap\OO_g(b)$ is finite.
\end{theorem}

\noindent
Hypothesis~(1) ensures that $f$ and $g$ do not share a common iterate. Obviously, if $f$ and $g$ share a common iterate then $\OO_f(a)\cap\OO_g(b)$ may be infinite (for example, when $a=b$). Hypothesis (2) means there does not exist any invertible $\mu \in k(z)$ for which both $\mu \circ f \circ \mu^{-1} \in K(z)$ and $\mu(a) \in \P^1(K)$.  We believe Theorem~\ref{main result application} should hold without hypothesis (2); we use it to guarantee that if $a$ is not preperiodic for $f$, then $\hhat_f(a)>0$ (\cite{Baker:functionfields}, \cite{D:stableheight}).  As mentioned above, we suspect that hypothesis (3) is always satisfied; the rationality of heights provides the decisive Diophantine ingredient needed in our proof.  

Finally, we note that essentially the same proof as the one we employ for Theorem~\ref{main result application} would work under the more general hypotheses that $f$ and $g$ are polarizable endomorphisms of degrees $d$ and $e$ with respect to the same ample line bundle on a projective variety $Z$.  Condition~(2) would be replaced with the assumption that $f$ does not descend to an endomorphism defined over the constant field $K$.  Then one would use \cite{CH1} and \cite{CH2} instead of \cite{D:stableheight} or \cite{Baker:functionfields} to infer that $\hhat_f(a)>0$ if $a$ is not preperiodic under the action of $f$. However, since checking whether hypothesis~(3) is always met for rational functions appears to be difficult, one might expect this same condition for the general case of polarizable endomorphisms to be very challenging.

\subsection{Acknowledgements}
We are grateful to Jason Bell, Rob Benedetto, Greg Call, Holly Krieger, and Joe Silverman for useful conversations regarding this project.  We also thank the anonymous referees of an earlier version of this manuscript for their helpful comments.

\bigskip
\section{Notation}

We fix some notation and terminology that we use throughout the article. 

\subsection{Places of function fields}

Let $K$ be an algebraically closed field of characteristic $0$, and let $k$ be the function field of a projective, smooth curve $X$ defined over $K$. We let $M_k$ be the set of all \emph{places} of $k$ which are trivial on $K$; each place  $v\in M_k$ corresponds to a point $P(v)\in X(K)$. By abuse of notation, we identify each place $v\in M_k$ with both its valuation function $v(\cdot )$ and also with its corresponding nonarchimedean absolute value $|\cdot |_v$ (note that $|x|_v:=e^{-v(x)}$ for each $x\in k$).  We  normalize each valuation $v$  so that $v(k)=\Z\cup\{\infty\}$ (note that $v(0)=\infty$). In other words, for a nonzero $x\in k$, the valuation $v(x)$ represents the order of vanishing at $P(v)$ of the function $x$ (seen as a map $X\mapsto \bP^1$).  
We have the following \emph{product formula} (see \cite{BG}) for each nonzero $x\in k$: 
\begin{equation}
\label{equation product formula}
\prod_{v\in M_k}|x|_v=1.
\end{equation}
Note that in \eqref{equation product formula}, each local degree (usually denoted by $N_v$) equals $1$ since the residue field at each place $v$ is canonically isomorphic to $K$. 

For each place $v\in M_k$ we  define $r_v:\bP^1(k)\lra \bP^1(K)$ to be the reduction map modulo $v$. So, for each $P\in \bP^1(k)$ we let $a,b\in k$ such that $\min\{v(a), v(b)\}=0$ and $P=[a:b]$; then $r_v(P)$ is obtained by reducing modulo $v$ both $a$ and $b$. Hence, for $x\in k$ such that $|x|_v\le 1$, $r_v(x)$ is simply the image of $x$ in the residue field at the place $v$, while if $|x|_v>1$, then $r_v(x)=\infty$.

For each $v\in M_k$ we let $k_v$ be the completion of $k$ at the place $v$; we let $\C_v$ be the completion of a given algebraic closure of $k_v$. Then $\C_v$ is a complete, algebraically closed field, and we fix an embedding of $\bar{k}$ into $\C_v$. We extend naturally the reduction map $r_v:\bP^1(\C_v)\lra \bP^1(K)$.

We let $\O_v = \{x\in \C_v: |x|_v \leq 1\} = \{x: v(x) \geq 0\}$.  
For each $a\in \C_v$ and each $r\ge 0$, we let 
	$$D(a,r) = \{x \in \C_v: |x-a|_v \leq r\}$$
and we refer to a set of the form $\{x:  |x-a| < r\}$ as an {\em open disk}.  

\subsection{Weil height}

We define the Weil logarithmic height function $h(\cdot )$ on $\P^1(\kbar)$ as follows: for each $P:=[x_0:x_1]\in \P^1(\kbar)$ we let 
$$h(P):=\frac{1}{[k(x_0,x_1):k]}\; \sum_{v\in M_k} \;  \sum_{\substack{\sigma:k(x_0,x_1)\hookrightarrow \kbar\\ \sigma|_k ={\rm id}|_k}} \; \log\max\left\{|\sigma(x_0)|_v, |\sigma(x_1)|_v\right\}.$$
We refer the reader to \cite{BG} for more details about the Weil height for function fields.

\subsection{Rational functions on function fields}  \label{ff language}
Let $f\in k(z)$ be a rational function of degree $d\ge 2$. We consider homogeneous presentations 
	$$F(z,w) = (P(z,w), Q(z,w))$$
of the rational function $f(z)= P(z,1)/Q(z,1)$, where $P$ and $Q$ are homogeneous polynomials of degree $d$.  The presentation $F$ is said to be {\bf normalized (at the place $v$)}  if the coefficients of $P$ and $Q$ lie in $\O_v$, and if at least one coefficient of $P$ and of $Q$ has valuation $0$.  The rational function determined by $F \mod v$ is called the {\bf reduction} of $f$; it has degree between $0$ and $d = \deg f$.  The map $f$ is said to have {\bf good reduction} if the degree is equal to $d$.  The map $f$ is said to have {\bf potential good reduction} if, after a change of coordinates by conjugation, the resulting map has good reduction.  See also \cite{Silverman:dynamics} for more information on these fundamental concepts.  

\subsection{Canonical height}  
\label{properties}

The canonical height of the dynamical system $f: \P^1\to \P^1$ defined over $k$ of degree $d\geq 2$ is defined by 
	$$\hhat_f(P) = \lim_{n\to\infty} \frac{1}{d^n} h(f^n(P)).$$
There exists a local decomposition associated to a choice of divisor on $\P^1$; the standard choice is to take $\infty$ as a divisor.  That is, if we write $f(z) = P(z)/Q(z)$, where $P$ and $Q$ have no common roots, then there is a unique function $\lhat_{f,v}$ satisfying \cite[Theorem 3.27]{Silverman:dynamics}:
\begin{enumerate}
\item  for all $a \in \P^1(\C_v) \setminus\{\infty\}$ with $Q(a) \not=0$, 
	$$\lhat_{f,v}(f(a)) = d \, \lhat_{f,v}(a) - \log|Q(a)|_v,$$
	and
\item $a \mapsto \lhat_{f,v}(a) - \log \max\{|a|_v, 1\}$ extends to a bounded continuous function on $\P^1(\C_v)$.   
\end{enumerate}
Note that if we replace the pair $(P,Q)$ with $(sP,sQ)$ for $s \in k^*$, the local height is transformed by adding 
	$$\frac{1}{d-1} \log|s|_v = -\, \frac{1}{d-1} \, v(s) \; \in \; \Q$$ 
\cite[Remark~3.28]{Silverman:dynamics}.  

For a polynomial $f\in k[z]$ of degree $d\ge 2$ and for any $v\in M_k$, the local canonical height $\lhat_{f,v}$ is given by
$$\lhat_{f,v}(x) =-\lim_{n\to\infty}\frac{\min\{0,v(f^n(x))\}}{d^n}$$
upon taking denominator $Q = 1$.  

As shown in \cite{Call:Silverman}, if $x = (x_0:x_1) \in \bP^1(\kbar)$ then 
	$$\hhat_f(x)=\frac{1}{[k(x_0,x_1):k]}\;  \sum_{v\in M_k}  \; \sum_{\substack{\sigma:k(x_0,x_1)\hookrightarrow \kbar\\ \sigma|_k ={\rm id}|_k}} \; \lhat_{f,v}(x).$$ 
Also, as proved in \cite{Call:Silverman}, the difference $|\hhat_f(x)-h(x)|$ is uniformly bounded over $x\in \P^1(\kbar)$ (this property will be used in Section~\ref{section diophantine}). Finally, we note the result of Baker \cite{Baker:functionfields}, who proved that unless $f$ is \emph{isotrivial} (i.e., there exists a linear map $\mu\in \kbar(z)$ such that $\mu^{-1}\circ f\circ \mu\in K(z)$), we have that $\hhat_f(x)=0$ if and only if $x$ is preperiodic (this property will be used in Section~\ref{section diophantine}).

\bigskip
\section{Computing the local canonical height}

We let $f\in k(z)$ be a rational function of degree $d\ge 2$.  In this section, we compute the local canonical height $\lhat_{f,v}$ following the approach of \cite[Section 5]{Call:Silverman} and \cite{Call:Goldstine}, and we explain how to determine its rationality.  We also show rationality does not depend on choices of coordinates, so in particular, if $f$ has potential good reduction, then all points have rational canonical height.

\subsection{Orders of vanishing} 
Take a homogeneous presentation of $f$, writing 
	$$F(z,w) = (P(z,w), Q(z,w))$$
as in \S\ref{ff language}.  Fix $v \in M_k$.  Given any point $a \in \P^1(\C_v)$, we express $a$ in homogeneous coordinates $A = (A_1, A_2) \in (\C_v)^2$ so that $\min\{v(A_1), v(A_2)\} = 0$.  We define {\bf orders} 
\begin{equation} \label{order}
	\sigma(F, a) = \min\{v(P(A)), v(Q(A)) \}.
\end{equation}
Note this is well defined, as it does not depend on our choice of representative of $A$, subject to the chosen normalization condition.  It does depend on the choice of $F$, since 
	$$\sigma(\alpha F, a) = \sigma(F,a) + v( \alpha)$$
for any $\alpha\in k^*$.

The function $\sigma(F, \cdot)$ is uniformly bounded on $\P^1(\C_v)$, as first observed in \cite{Call:Silverman}; in fact, for normalized $F$, we have
\begin{equation} \label{order bound}
	0 \leq \sigma(F,a) \leq v(\Res(F))
\end{equation}
for all $a \in \P^1(\C_v)$ \cite[Lemma 10.1]{BRbook},  
where $\Res(F)$ is the homogeneous resultant of the pair $(P,Q)$. In particular, if $f$ has good reduction at $v$, which is equivalent with $v(\Res(F))=0$ then $\sigma(F,a)=0$ for all $a$ in $\P^1(\C_v)$.  Also, restricting to $k$ or its completion $k_v$ (or any finite extension), the function $\sigma(F, \cdot)$ takes only finitely many rational values.  

From the definitions, we compute that 
	$$\sigma(F^2, a) = d \, \sigma(F,a) + \sigma(F, f(a)),$$
and 
\begin{eqnarray}  \label{order iterate}
\sigma(F^n, a) &=& d \, \sigma(F^{n-1}, a) + \sigma(F, f^{n-1}(a)) \\
	&=& d^{n-1}\, \sigma(F, a) + d^{n-2} \, \sigma(F, f(a)) + \cdots + \sigma(F, f^{n-1}(a)) \nonumber
\end{eqnarray}
for all $n\geq 1$.  The uniform bound on $\sigma(F,\cdot)$ shows that the limit 
\begin{equation} \label{eta}
	\eta(F,a) := \lim_{n\to\infty}  \frac{\sigma(F^n, a)}{d^n} = \frac{1}{d} \sum_{n=0}^\infty \frac{\sigma(F, f^n(a))}{d^n}
\end{equation}
exists.  

\begin{lemma} {\em (\cite[Theorem 5.3]{Call:Silverman}, \cite[\S3]{Call:Goldstine})} \label{local height}
The local canonical height (relative to $\infty$) for $F = (P,Q)$ is computed as 
	$$\lhat_{f,v}(a) = -\eta(F,a) - \min\{0, v(a)\}$$
for all $a \not= \infty$.
\end{lemma}

The rationality of $\lhat_{f,v}(a)$ thus hinges on the values of the sequence
\begin{equation}
\label{boundedness of the c_n}
\sigma_n := \sigma(F, f^n(a)) \; \in \; \Q \cap [0, v(\Res F)].
\end{equation}
We note the following basic fact about $d$-ary expansions for the proofs of our main theorems: 

\begin{prop}  \label{rationality test}
Fix integer $d\geq 2$.  Suppose that $\sigma_n \in \{0, 1, \ldots, d-1\}$ for all $n$.  Then 
	$$\eta = \frac{1}{d} \sum_{n=0}^\infty \frac{\sigma_n}{d^n}$$
is rational if and only if the sequence $\{\sigma_n\}$ is eventually periodic.
\end{prop}

\subsection{Transforming the local canonical height}
\label{transforming}
Recall from \S\ref{properties} that $\lhat_{f,v}$ depends on the choice of polynomials $P$ and $Q$ representing $f = P/Q$.  However, regardless of the choice, the property whether $\lhat_{f,v}(\cdot )$ takes rational values or not is unchanged. 

Next we show that replacing a rational function $\Phi(z) = P(z)/Q(z)$ by conjugate rational function $\Psi$, i.e., for which there exists a M\"obius transformation $\eta\in k(z)$ such that $\Phi=\eta^{-1}\circ \Psi\circ \eta$, the property of taking only rational values by the local canonical height is unaffected. Since each linear transformation $\eta$ is a composition of finitely many linear transformations of the form
\begin{enumerate}
\item[(1)] $z\mapsto cz$ for some constant $c\in k^*$;
\item[(2)] $z\mapsto z+c$ for some constant $c\in k^*$; and
\item[(3)] $z\mapsto 1/z$,
\end{enumerate}
it suffices to consider separately each of the cases (1)--(3).

\smallskip\noindent
{\bf Case (1):} The change of coordinates given by multiplicative translation.  For each place $v$, we have
that the local canonical heights are related by
\begin{equation}
\label{formula 1}
\lhat_{\Phi, v}(x) = \lhat_{\Psi, v}(cx) + v(c).  
\end{equation}
The proof goes through checking the properties (1) and (2) of \S\ref{properties} for the local canonical height. Since $\Psi(z) = P(z)/Q(z)$, then $\Phi(z) = P(cz)/(c\cdot Q(cz))$. 
We take the denominator of $\Phi(z)$ to be $c^{1-d}Q(cz)$.  Using the expression of \eqref{formula 1} as a definition for $\lhat_{\Phi,v}(\cdot)$ in terms of $\lhat_{\Psi,v}(\cdot)$, we see that
\begin{eqnarray*}
\lhat_{\Phi,v}(x) + \min\{v(x), 0\}
& =& \lhat_{\Psi, v}(cx) + v(c) + \min\{v(x), 0\}\\        
& = &- \min\{v(cx), 0\} + O(1) + \min\{v(cx), v(c)\}\\  
& = & O(1),
\end{eqnarray*}
where in the next to the last equality  we used property~(2) for $\lhat_{\Psi,v}$. This shows $\lhat_{\Phi,v}$ satisfies property (2) as well. Now, for property (1) for $\lhat_{\Phi,v}$, we have
\begin{eqnarray*}
\lhat_{\Phi,v}(\Phi(x))
& = & \lhat_{\Psi, v}(c\Phi(x)) + v(c)\\    
& = & \lhat_{\Psi, v}(\Psi(cx)) + v(c)\\
& =  & d \lhat_{\Psi,v}(cx) - \log|Q(cx)|_v + v(c)\\   
& = & d \lhat_{\Phi, v}(x) - \log|c\cdot Q(cx)|_v,
\end{eqnarray*}
which yields formula (ii) for $\lambda_{\Phi, v}$ since $-\log|c|_v=v(c)$ and the denominator for $\Phi(z)$ is $c^{1-d}Q(cz)$.

\smallskip\noindent
{\bf Case (2):} The change of coordinates given by additive translation.  In this case, 
$$\Phi(z)=\Psi(z+c)-c=\frac{P(z+c)}{Q(z+c)}-c=\frac{P(z+c)-cQ(z+c)}{Q(z+c)}.$$
We observe that the denominator $Q(z+c)$ remains coprime with the numerator $P(z+c)-cQ(z+c)$; also, $\max\{\deg(P(z+c)-cQ(z+c)), \deg(Q(z+c))\}=d$. We claim that for any place $v$, we have 
\begin{equation} \label{formula 2}
\lhat_{\Phi,v}(x)=\lhat_{\Psi,v}(x+c).
\end{equation}
Again, we only need to check that $\lhat_{\Phi,v}$ satisfies the two conditions~(1)--(2) above. Condition~(2) follows because
\begin{eqnarray*}
\lhat_{\Phi,v}(x)+\min\{v(x),0\}
& =& \lhat_{\Psi,v}(x+c)+\min\{v(x),0\}\\
& =& \left(\lhat_{\Psi,v}(x+c)+\min\{v(x+c),0\}\right) + \left(\min\{v(x),0\} - \min\{v(x+c), 0\}\right)\\
& =& O(1).
\end{eqnarray*}
The first parenthesis above is $O(1)$ using condition~(2) applied to the local canonical height $\lhat_{\Psi,v}$.  
Now, for condition~(2), we check that
\begin{eqnarray*}
\lhat_{\Phi,v}(\Phi(x))
& =& \lhat_{\Psi,v}(\Phi(x)+c)\\
& =& \lhat_{\Psi,v}(\Psi(x+c))\\
& =& d\cdot \lhat_{\Psi,v}(x+c) - \log|Q(x+c)|_v,
\end{eqnarray*}
which proves condition~(1) for $\lhat_{\Phi,v}$ because the denominator of $\Phi(z)$ is indeed $Q(z+c)$.

\smallskip\noindent
{\bf Case (3):} The change of coordinates given by inversion.  Since $\Phi(z)=1/\Psi(1/z)$, then 
$$\Phi(z)=\frac{Q\left(\frac{1}{z}\right)}{P\left(\frac{1}{z}\right)}= \frac{z^d\cdot Q\left(\frac{1}{z}\right)}{z^d\cdot P\left(\frac{1}{z}\right)}.$$
So, the denominator of $\Phi(z)$ is $z^d\cdot P(1/z)$, which is a polynomial coprime with the  polynomial $z^d\cdot Q(1/z)$, which is the numerator of $\Phi(z)$. We claim that for any place $v$, we have that 
\begin{equation}
\label{formula 3}
\lhat_{\Phi,v}(x)=\lhat_{\Psi,v}\left(\frac{1}{x}\right) + v\left(\frac{1}{x}\right).
\end{equation}
Now, to check condition~(2) for the above defined function $\lhat_{\Phi,v}(\cdot )$, we compute
\begin{eqnarray*}
\lhat_{\Phi,v}(x)+\min\{v(x),0\}
& =& \lhat_{\Psi,v}\left(\frac{1}{x}\right) + v\left(\frac{1}{x}\right) + \min\{v(x), 0\}\\
& =& \lhat_{\Psi,v}\left(\frac{1}{x}\right) + \min\left\{0, v\left(\frac{1}{x}\right)\right\}\\
& =& O(1),
\end{eqnarray*}
because $\lhat_{\Psi,v}$ satisfies condition~(2). For condition~(1), we have that
\begin{eqnarray*}
\lhat_{\Phi,v}(\Phi(x))
& =& \lhat_{\Psi,v}\left(\frac{1}{\Phi(x)}\right) + v\left(\frac{1}{\Phi(x)}\right)\\
& = & \lhat_{\Psi,v}\left(\Psi\left(\frac{1}{x}\right)\right) + v\left(\Psi\left(\frac{1}{x}\right)\right)\\
& = & d\cdot \lhat_{\Psi,v}\left(\frac{1}{x}\right) - \log\left|Q\left(\frac{1}{x}\right)\right|_v + v\left(P\left(\frac{1}{x}\right)\right) - v\left(Q\left(\frac{1}{x}\right)\right)\\
& =& d\cdot \left(\lhat_{\Psi,v}\left(\frac{1}{x}\right) + v\left(\frac{1}{x}\right)\right) +v\left(x^d\right) + v\left(P\left(\frac{1}{x}\right)\right)\\
& =& d\cdot \lhat_{\Phi,v}(x) + v\left(x^d\cdot P\left(\frac{1}{x}\right)\right)\\
& =& d\cdot \lhat_{\Phi,v}(x) - \log\left|x^d\cdot P\left(\frac{1}{x}\right)\right|_v,
\end{eqnarray*}
as desired.

\begin{prop} \label{potential good reduction rationality}
If rational function $f \in k(z)$ of degree $\geq 2$ has potential good reduction at the place $v$, then $\lhat_{f,v}(a) \in \Q$ for all $a \in \P^1(\C_v)\setminus\{\infty\}$.
\end{prop}

\proof
Let $F$ be a normalized homogeneous presentation of $f$.  Suppose that $f$ has good reduction.  As observed immediately after \eqref{order bound}, good reduction implies that $\eta(F,a) \equiv 0$.  Therefore, from Lemma \ref{local height}, we have $\lhat_{f,v}(a)=-\min\{0, v(a)\}\in \Q$ for all $a \in\P^1(\C_v)\setminus \{\infty\}$.  On the other hand, if $f$ has potential good reduction, then equations \eqref{formula 1}, \eqref{formula 2}, and \eqref{formula 3} show that $\lhat_{f,v}$ will again take only rational values.
\qed

\bigskip
\section{The spine in Berkovich space}

We begin this section with the background we need and the terminology we use about the Berkovich projective line $\bP^1_v$ over $\C_v$.  We introduce the {\bf spine} of a rational function $f$ inside $\bP_v^1$, a finite tree that captures the geometry of its zeroes and poles.  It is our tool for computing the order function \eqref{order} for arbitrary rational functions $f$.  

\subsection{Berkovich space and definition of the Julia set}  \label{Berkovich background}
Useful introductions to the Berkovich projective line $\bP^1_v$ are given in \cite{BRbook} and \cite{Faber:ramification1}.  See \cite[Sections 1 and 2]{Kiwi:quad} for a presentation of the relevant dynamical notions.  We summarize the most important things here.  

In the weak topology, the projective line $\P^1(\C_v)$ is a dense subset of $\bP^1_v$.  Its points are called Type I or classical.  The closure inside $\bP^1_v$ of a disk $D(a,r) = \{z: |z-a|_v \leq r\}$ in $\C_v$ is called a Berkovich closed disk.  It has a unique boundary point inside $\bP^1_v$.  If the radius $r$ lies in the value group $|\C_v^*| = \Q$, this boundary point is called a Type II point, and it is a branching point of the tree structure on $\bP^1_v\setminus\P^1(\C_v)$.  In this way, we identify Type II points with disks of rational radius; we sometimes write $\zeta = D(a,r)$ where $\zeta$ is the Type II point bounding the disk $D(a,r)$.  Given a Type II point $\zeta \in \bP^1_v$, the connected components of $\bP^1_v\setminus\{\zeta\}$ are called {\bf directions from} $\zeta$.   The intersection of a direction with the set of Type I points is an open disk. There is a distinguished Type II point $\zeta_0$ called the {\bf Gauss point} corresponding to the unit disk $D(0,1)$.  The set of directions from $\zeta_0$ is isomorphic to $\P^1(\O_v/\mathfrak{m}_v) = \P^1(K)$.

The Type II points may also be identified with M\"obius changes of coordinates on $\P^1(\C_v)$.  More precisely, the group $\PSL_2(\C_v)$ acts on $\bP^1_v$, and the stabilizer of the Gauss point $\zeta_0$ is equal to $\PSL_2(\O_v)$.  For any Type II point $\zeta$ of $\bP^1_v$, there is a unique $M \in \PSL_2(\C_v)$, modulo postcomposition by an element of the stabilizer, that sends $\zeta$ to $\zeta_0$.  See the discussion in \cite[\S3]{Kiwi:rescaling}.

Now suppose that $f$ is a rational function of degree $d\geq 2$ defined over $k$.  As a map, $f$ takes (sufficiently small) disks in $\P^1(\C_v)$ to disks; the image of an arbitrary disk will be either another disk or all of $\P^1(\C_v)$.  In this way $f$ extends naturally and continuously to $\bP^1_v$.  Topologically, it behaves much like a holomorphic branched cover of the Riemann sphere; see \cite{Faber:ramification1, Faber:ramification2} for a study of the branching structure.  See also \cite[Proposition 2.1]{DF:degenerations2} for a description of its local topological degrees, exactly analogous to the behavior of proper analytic maps on regions in the complex plane, built on a non-archimedean version of the Argument Principle.   

Given any Type II point $\zeta$, the action of $f$ on its set of directions is easily computed.  Indeed, under a pair of M\"obius changes of coordinates $A$ and $B$, we can move $\zeta$ and its image $f(\zeta)$ to the Gauss point, so that $B \circ f \circ A^{-1}$ fixes the Gauss point $\zeta_0$.  The reduction modulo $v$ of $B \circ f \circ A^{-1}$ (expressed in a homogeneous, normalized form) determines a rational function of degree $d(\zeta) \geq 1$ which determines the action from the set of directions at $\zeta$ mapping to those of $f(\zeta)$.   In particular, the set of directions from $\zeta$ always maps surjectively over the set of directions from $f(\zeta)$.  Moreover, a direction will always map surjectively over its image direction; but there are particular directions from $\zeta$ that will map over all of $\bP^1_v$.  These are called {\bf bad} directions, and they arise when the degree of the map defined by $B \circ f \circ A^{-1} \;\rm{mod} \; v$ satisfies $d(\zeta) < d$; a bad direction at $\zeta$ is sent by $A$ to a direction from $\zeta_0$ containing both a pole and a zero of $B \circ f \circ A^{-1}$. 

A point $x \in \bP^1_v$ belongs to the {\bf Julia set} of $f$ if, for every open neighborhood $V$ containing $x$, the union $\bigcup_{n\geq 0} f^n(V)$ omits at most two points of $\bP^1_v$.  Otherwise, $x$ is said to be in the {\bf Fatou set} of $f$.  For example, a map $f$ has good reduction if and only if its Julia set is equal to the Gauss point $\{\zeta_0\}$.  It has potential good reduction if and only if the Julia set consists of a single Type II point.  


\subsection{The spine of a rational map}  
Given a rational function $f$ defined over $k$, write $f$ in homogeneous coordinates as 
	$$F(z,w) = (P(z,w), Q(z,w)).$$
Recall that $F$ is normalized if all coefficients lie in $\O_v$ and if $P$ and $Q$ each have at least one coefficient of valuation 0.  Recall the order function $\sigma(F, \cdot)$ on $\P^1(\C_v)$, defined above in \eqref{order}:  
	$$\sigma(F, a) = \min\{v(P(A)), v(Q(A)) \}$$
where $A = (A_1, A_2)$ is any expression for $a\in \P^1(\C_v)$ in homogeneous coordinates with $\min\{v(A_1), v(A_2)\} = 0$.  The following result, proved in \S\ref{order function proof}, allows us to define the spine of a rational function.

\begin{prop} \label{order function}
There is a well-defined and continuous function 
	$$\zeta \mapsto \sigma_F(\zeta)$$ 
on Type II points $\zeta$ in $\bP^1_v$, satisfying: 
\begin{itemize}
\item[(i)]	$\sigma_F(\zeta_0) = 0$ for the Gauss point $\zeta_0$ if and only if $F$ is normalized;
\item[(ii)]	the function $\sigma_F$ is monotone (non-strictly) increasing along any path from $\zeta_0$ to a zero or pole of $f$ and locally constant elsewhere; and 
\item[(iii)] for all but finitely many directions $U$ from $\zeta$, and every point $a \in U\cap \P^1(\C_v)$, the order $a \mapsto \sigma(F,a)$ is constant equal to $\sigma_F(\zeta)$.
\end{itemize}
\end{prop}

The {\bf spine} $S_f$ is the smallest, closed, connected subset of $\bP^1_v$, containing the Gauss point $\zeta_0$, so that the order function $\sigma_F$ is constant on connected components of $\bP^1_v\setminus S_f$.  From Proposition \ref{order function}, the spine is contained in the finitely-branched tree spanning the zeroes and poles of $f$ and containing the Gauss point.  It is not a dynamical object, and it is not a conjugacy invariant.  We will show:

\begin{theorem}  \label{spine shape}
The spine $S_f$ is the connected hull of the Gauss point $\zeta_0$ and its set of preimages $f^{-1}(\zeta_0)$ in $\bP^1_v$.  
\end{theorem}

\noindent
In particular, the spine $S_f$ is equal to the Gauss point $\{\zeta_0\}$ if and only if $f$ has good reduction.

\subsection{Computing the order function}  \label{order function proof}
Fix a normalized homogeneous presentation of $F$ as 
\begin{equation} \label{general F}
F(z,w) = (P(z,w), Q(z,w)) = \left(c \prod_{i=1}^d (\beta_i z - \alpha_i w), u \prod_{j = 1}^{d} (\delta_j z - \gamma_j w) \right)
\end{equation}
with $\min\{v(c), v(u)\} = \min\{v(\alpha_i), v(\beta_i)\} = \min\{v(\delta_j), v(\gamma_j)\} = 0$ for all $i$ and $j$.  Observe that $\{z_i = \alpha_i/\beta_i\}$ is the set of zeroes of $f$ and $\{p_j = \gamma_j/\delta_j\}$ is the set of poles of $f$.  

Each point $a$ in the disk $D(0,1)$ in $\P^1(\C_v)$ will be represented in homogeneous coordinates as $A = (a, 1)$.  Then 
	$$\sigma(F, a) := \min\{v(P(A)), v(Q(A))\}$$ 
can be computed by observing that 
	$$v(\beta_i x - \alpha_i) = \left\{  \begin{array}{ll} 
				0 & \mbox{if } |z_i|_v > 1 \\
				v(a - z_i) & \mbox{if } |z_i|_v \leq 1  \end{array} \right. $$
and 
	$$v(\delta_i x - \gamma_i) = \left\{  \begin{array}{ll} 
				0 & \mbox{if } |p_i|_v > 1 \\
				v(a - p_i) & \mbox{if } |p_i|_v \leq 1  \end{array} \right. $$
Therefore, still assuming that $|a|_v \leq 1$, we have
\begin{equation} \label{small a}
\sigma(F, a) = \min\left\{ v(c) + \sum_{ \{i: |z_i|_v \leq 1\} } v(a-z_i) \; ,\;  v(u) + \sum_{ \{i: |p_i|_v \leq 1\} } v(a-p_i) \right\}.  
\end{equation}

Similarly, a point $a$ with $|a|_v >  1$, allowing also $a=\infty$, may be represented as $A = (1, 1/a)$.  Then $\sigma(F, a)$ can be computed by observing that 
	$$v(\beta_i - \alpha_i y) = \left\{  \begin{array}{ll} 
				0 & \mbox{if } |z_i|_v \leq 1 \\
				v(1/a - 1/z_i) & \mbox{if } |z_i|_v > 1  \end{array} \right. $$
and 
	$$v(\delta_i  - \gamma_i y) = \left\{  \begin{array}{ll} 
				0 & \mbox{if } |p_i|_v \leq 1 \\
				v(1/a - 1/p_i) & \mbox{if } |p_i|_v > 1  \end{array} \right. $$
Consequently, for all $|a|_v > 1$, we compute that 
\begin{equation} \label{large a}
\sigma(F, a) = \min\left\{ v(c) + \sum_{ \{i: |z_i|_v > 1\} } v(1/a- 1/z_i) \; ,\;  v(u) + \sum_{ \{i: |p_i| > 1\} } v(1/a- 1/p_i) \right\}.  
\end{equation}

\medskip
\proof[Proof of Proposition \ref{order function}]
We will use condition (iii) of the Proposition to define $\sigma_F$ for the $F$ given by formula \eqref{general F}; it is a normalized homogeneous presentation for $f$.  Because $\min\{v(c), v(u)\} = 0$, the equations \eqref{small a} and \eqref{large a} show that $\sigma(F, a) =0$ unless $a$ lies in one of the finitely many directions from the Gauss point $\zeta_0$ which contains a pole or a zero.  Thus $\sigma_F(\zeta_0) = 0$.  Recalling that $\sigma(sF, a) = \sigma(F,a) + v(s)$ for any $s\in k^*$, we have the proof of (i).  

The equations \eqref{small a} and \eqref{large a} also show that the dependence of $\sigma(F,a)$ upon $a$ only incorporates its distance to the zeroes and poles of $f$.  In particular, points lying in all but finitely many directions from any Type II point $\zeta$ will have the same order.  This shows that $\sigma_F$ is well defined, satisfying (iii) of the Proposition, and it is continuous.  

Now fix any two points $x \not= y$ in the disk $D(0,1) = \{z: |z|_v \leq 1\}$.  The function 
	$$a \mapsto v(a - y)$$
is (non-strictly) increasing, for $a$ in $D(0,1)\setminus \{z: |z-y|_v < |x-y|_v\}$, as $v(a-x)$ increases from 0 to $\infty$.  Its maximum value over this domain is $v(x-y)$.  Therefore, for any Type II point of the form $\zeta = D(x,r)$ with $r \leq 1$, equation (\ref{small a}) shows that the order $\sigma_F$ can only increase (or remain constant) along the path from $\zeta_0$ towards $\zeta$.  Similarly for Type II points lying in the direction of $\infty$ from $\zeta_0$, using \eqref{large a}.  This proves (ii).  
\qed

\subsection{Proof of Theorem \ref{spine shape}}
We begin with a lemma about the preimages of the Gauss point.  Let $F = (P,Q)$ be a normalized presentation of $f$.   Recall that a point $a \in \P^1(\C_v)$ will be expressed in homogeneous coordinates as $A = (a,1)$ for $|a|_v \leq 1$ and $A = (1, 1/a)$ for $|a|_v > 1$.  

\begin{lemma} \label{Gauss preimage}
Let $\zeta_0$ be the Gauss point, and suppose $f(\zeta) = \zeta_0$.  Then for $a$ in all but finitely many directions from $\zeta$, we have $v(P(A)) = v(Q(A)) = \sigma_F(\zeta)$.  
\end{lemma}

\proof 
Recall that the space of directions from any Type II point $\zeta \in \bP^1_v$ is identified with $\P^1(K)$.  The action of $f$, from the space of directions from $\zeta$ to the space of directions from $\zeta_0$, is by a rational function of degree $\geq 1$.  Furthermore, the image of all but finitely many directions from $\zeta$ is equal to a direction from $\zeta_0$.  In particular, all but finitely many directions from $\zeta$ are sent to directions from $\zeta_0$ that do {\em not} contain 0 or $\infty$.  This implies that $|f(a)|_v = 1$ for any point $a$ in those directions, so for these $a$, we must have $v(P(A)) =v(Q(A))$.    
\qed

\medskip
Suppose that $v(c) = v(u) = 0$ in the expression \eqref{general F} for $F$.  Then, from the formulas \eqref{small a} and \eqref{large a}, we find that $\sigma(F, a) >0$ if and only if $a$ lies in a direction from $\zeta_0$ containing both a zero and a pole.  We now show that this condition is equivalent to there being a preimage of $\zeta_0$ in that direction.  Let $C$ denote the path joining a zero and pole lying in the same direction from $\zeta_0$.  By continuity, the image $f(C)$ must contain the path joining $0$ and $\infty$, so in particular it contains $\zeta_0$.  In other words, there is a preimage of $\zeta_0$ inside $U$.  Conversely, if a direction $U$ from $\zeta_0$ contains a preimage $\zeta$ of $\zeta_0$, then we examine the action of $f$ on directions from $\zeta$.  Each direction from $\zeta$ maps surjectively over its image direction (and possibly also over all of $\bP^1_v$), and the action on the space of directions is surjective.  Therefore, restricting to the directions from $\zeta$ that are contained in $U$ (which is all but one direction), we conclude that $U$ must contain at least one zero or pole of $f$.  So let us assume that the direction $U$ contains a zero $z_i$.  This implies that for all $a\in U$, we have $v(a - z_i) >0$ (or that $v(1/a - 1/z_i) > 0$ in the case where $U$ is the direction of $\infty$), so that $v(P(A))>0$.  But by Lemma \ref{Gauss preimage}, we know that $\sigma_F(\zeta) = v(Q(A)) = v(P(A)) = \sigma(F,a)$ for $a$ in all but finitely many directions from $\zeta$, and therefore $\sigma_F(\zeta)>0$.  Consequently, there must also be a pole lying in the direction $U$.  If $U$ contained a pole of $f$ rather than a zero, the argument is similar, so we see that $U$ contains both a pole and a zero.  

Let $T_f$ be the finite tree spanned by the Gauss point $\zeta_0$ and its finitely many preimages, with its vertices consisting of the elements of $\{\zeta_0\} \cup f^{-1}(\{\zeta_0\})$ together with any additional branching points in $T_f$.  We call a preimage $\zeta$ of $\zeta_0$ an {\bf end} of $T_f$ if $\zeta \not= \zeta_0$ and if the valence of $\zeta$ in $T_f$ is equal to 1.  The preceding paragraph shows that, under the assumption that $v(c) = v(u) = 0$, the tree $T_f$ is a subset of the span of the zeroes and poles and $\zeta_0$ inside $\bP^1_v$; also that $\sigma_F(\zeta)>0$ for each preimage $\zeta$ of $\zeta_0$.  It remains to show that $\sigma_F$ achieves its maximum at ends of $T_f$ and not before, when moving along a path from $\zeta_0$ towards an end.  

Note that there cannot be a zero and a pole in the same connected component of $\bP^1\setminus T_f$; if there were, then the continuity argument above shows the path between them contains another preimage of $\zeta_0$.  Let $\zeta$ be an end of $T_f$.  Recall that, for $a$ in all but finitely many directions from $\zeta$, we have that $v(P(A)) = v(Q(A))$ by Lemma \ref{Gauss preimage}.  If we take $a$ in any direction disjoint from $T_f$, then at most one of $v(P(A))$ or $v(Q(A))$ can be larger than the value of $\sigma_F(\zeta)$.  Therefore, the order function $\sigma_F$ must remain constant in that direction.  This shows that $\sigma_F$ achieves its maximum on $T_f$ and that the spine $S_f$ is a subset of $T_f$. 

On the other hand, suppose that $\zeta'$ is a Type II point lying on the path between $\zeta_0$ and an endpoint $\zeta$ of $T_f$.  We will see that $\sigma_F(\zeta') < \sigma_F(\zeta)$.  For $\zeta'$ sufficiently close to $\zeta$, the point $f(\zeta')$ lies in the direction from $\zeta_0$ determined by the reduction map from $\zeta$ to $\zeta_0$.  We consider three cases.  If this image direction does not point to either 0 or $\infty$, then necessarily the number of zeroes and poles beyond $\zeta$ (outside $T_f$) must coincide -- because the action on the space of directions by $f$ has a well-defined degree.  Thus, moving from $\zeta$ to $\zeta'$ will reduce the valuations $v(P(A))$ and $v(Q(A))$ simultaneously; this shows that the order function decreases.  If $f(\zeta')$ lies in the direction of 0 from $\zeta_0$, this means exactly that $v(P(A)) > v(Q(A))$ for points $a$ in the annulus between $\zeta$ and $\zeta'$.  On the other hand, there must be at least one pole in the directions from $\zeta$ disjoing from $T_f$, and therefore $v(Q(A))$ will have decreased from the value $\sigma_F(\zeta)$ for all $a$ in the annulus between $\zeta$ and $\zeta'$.  This shows that $\sigma_F(\zeta') < \sigma_F(\zeta)$.  Similarly, if $f(\zeta')$ lies in the direction of $\infty$ from $\zeta_0$.

Now suppose $v(c)>v(u)=0$ in the expression \eqref{general F} for $F$.  This means that $f(\zeta_0)$ lies in the direction of 0 from $\zeta_0$.  Then $\sigma(F,a)>0$ if and only if $a$ lies in a direction containing a pole, from the formulas \eqref{small a} and \eqref{large a} for $\sigma(F, \cdot)$.  But any direction containing a pole $p$ must also contain a preimage of $\zeta_0$, because the image of the path $[\zeta_0, p]$ must contain $\zeta_0$ by continuity.  Conversely, if a direction $U$ from $\zeta_0$ contains a preimage $\zeta$ of $\zeta_0$, it must also contain a pole, by the surjectivity of $f$ on directions.  Indeed, each direction from $\zeta$ maps surjectively over its image direction from $\zeta_0$.  The only concern is if the direction from $\zeta$ back to $\zeta_0$ is the unique direction sent towards $\infty$.  But since $\zeta_0$ itself is moving towards 0, it would follow that the direction is a bad direction (otherwise its image would be exactly equal to the direction from $\zeta_0$ to infinity) -- so that direction maps surjectively over all of $\bP^1_v$ and by surjectivity of directions at all Type II points between $\zeta$ and $\zeta_0$, it follows that there must be a pole between $\zeta$ and $\zeta_0$.  

Similarly, if $0 = v(c) < v(u)$, then $f(\zeta_0)$ lies in the direction of $\infty$ from $\zeta_0$.  It follows that $\sigma(F,a)>0$ if and only if $a$ lies in a direction containing a zero, and again this is equivalent to the statement that the direction contains a preimage of $\zeta_0$.  
\qed

\subsection{Spine examples}
We conclude this section by computing the spine in a few examples.  We do so directly, with the formulas \eqref{small a} and \eqref{large a}, using the definition of the spine.  We also do the computation indirectly, appealing to Theorem \ref{spine shape}.  The first example has potential good reduction.  The second is a quadratic polynomial that does not have potential good reduction.

\begin{example}
Let $k = K(t)$ and $v = \ord_{t=0}$.  Consider 
	$$g(z) = \frac{z^2-1}{z}.$$
This function $g$ has good reduction, so its spine $S_g$ is the singleton $\{\zeta_0\}$ and the order function is $\equiv 0$.  Conjugating by $A(z) = t z$, we define
	$$f(z) = A \circ g \circ A^{-1} (z) = \frac{z^2 - t^2}{z}.$$
The zeroes of $f$ lie at $z_1 = t$ and $z_2 = -t$.  The poles of $f$ lie at 0 and $\infty$.  In the form \eqref{general F}, a homogeneous presentation of $f$ is given by
	$$F(z,w) = (z^2 - t^2 w^2, z w) = ((z-tw)(z+tw), z w)$$
with $c = u = 1$.  It is immediate from \eqref{large a} that $\sigma(F,a) = 0$ for all $|a|_v >1$.  Using \eqref{small a}, we compute that 
	$$\sigma(F, a) = \min\{v(a-t) + v(a+t), \; v(a)\}$$
for each $a \in \P^1(\C_v)$ with $|a|_v \leq 1$.  Consequently, 
	$$\sigma(F,a) = \left\{ \begin{array}{ll}  0 & \mbox{for } v(a) \leq 0 \\
								v(a) & \mbox{for } 0 \leq v(a) \leq 2 \\
								2  & \mbox{for } 2 \leq v(a) \end{array}  \right.$$
and the spine $S_f$ is the interval $[\zeta_0, \zeta_2]$ in $\bP_v^1$, from the Gauss point $\zeta_0$ and the point $\zeta_2$ representing $D(0, |t^2|_v)$.  

To compute the spine of $f$ with Theorem \ref{spine shape}, we need only determine the preimages of the Gauss point $\zeta_0$ under $f$.  Observing first that 
	$$f(z) = z + O(t^2)$$
we see that $f$ must fix the Gauss point; the reduction map has degree 1.  We also compute that 
	$$f(t^2 z) = -1/z + O(t^2),$$
which shows that the point $\zeta_2 = D(0, |t^2|_v)$ is another preimage of $\zeta_0$.  Because $f$ has degree 2, there are the only preimages of $\zeta_0$.  Thus, the spine $S_f$ is the convex hull of $\zeta_0$ and $\zeta_2$.  
\end{example}

\begin{example}
Let $k = K(t)$ and $v= \ord_{t=0}$.  Consider 
	$$f(z) = t^{-3} z(z-t).$$
In the form \eqref{general F}, a homogeneous presentation of $f$ is given by 
	$$F(z,w) = (z(z - t w), t^3 w^2)$$
so that $c = 1$ and $u = t^3$ has $v(u) = 3$.  For $|a|_v> 1$, we compute from \eqref{large a} that $\sigma(F, a) = 0$, as the zeroes of $f$ lie at 0 and $t$.  For $|a|_v \leq 1$, we have
	$$\sigma(F, a) = \min\{v(a) + v(a-t), \; 3\}.$$
We deduce that the spine $S_f$ has the shape of a Y.   The three tips of the Y lie at the Gauss point $\zeta_0$, at the point $\zeta_{0,2}$ representing the disk $D(0, |t^2|_v)$, and at the point $\zeta_{t,2}$ representing the disk $D(t, |t^2|_v)$, with the center of the Y at $\zeta_{0,1}$ representing the disk $D(0, |t|_v)$.  The function $\sigma_F$ of Proposition \ref{order function} satisfies $\sigma_F(\zeta_0) = 0$,  $\sigma_F(\zeta_{0,1}) = 2$, $\sigma_F(\zeta_{0,2}) = 3$, and $\sigma_F(\zeta_{t,2}) = 3$.  The function $\sigma_F$ is linear along the three branches joining the end points to the center $\zeta_{0,1}$ of the Y.

To compute the spine $S_f$ using Theorem \ref{spine shape}, we observe that 
	$$f(t^2 z) = -z + O(t)$$
and 
	$$f(t + t^2 z) = z + O(t).$$
These computations show that the points $\zeta_{0,2} = D(0, |t^2|_v)$ and $\zeta_{t,2} = D(t, |t|^2)$ are both sent by $f$ to $\zeta_0$.  As $f$ has degree $2$, these are the only preimages of $\zeta_0$.  Therefore, the spine is the convex hull of the set $\{\zeta_0, \zeta_{0,2}, \zeta_{t,2}\}$.  
\end{example}

\bigskip
\section{Elliptic curves, dynamically}

In this section we present our dynamical proof of Theorem \ref{dynamical E}.

\subsection{The Latt\`es tent}
\label{subsection Lattes}

Let $E$ be an elliptic curve over function field $k = K(X)$.  Via the quotient identifying a point $P$ with its inverse $-P$, the multiplication-by-2 endomorphism on $E$ induces a rational function $f: \P^1\to \P^1$ defined over $k$ of degree 4.  The N\'eron-Tate canonical height on $E$ is equal to $\frac12$ times the dynamical canonical height of $f$ (see \eqref{definition canonical height elliptic curve} and Subsection~\ref{properties}).  

For a given place $v$, either $E$ has potential good reduction, or the $j$-invariant satisfies $|j|_v>1$ and we may present $E$ locally as $\C_v^*/q^\Z$ for some $|q|_v < 1$ \cite[Chapter V Theorem 5.3(a), Remark 3.1.2]{Silverman:Advanced}.  In the latter case, from \cite[Proposition 5.1]{FRL:ergodic} and its proof, we know that the action of $f$ on the Berkovich projective line $\bP^1_v$ has a very simple form. 

Indeed, after a change of coordinates (putting our $E$ into Legendre form), we may assume that our $f$ is 
	$$f(z) = \frac{(z^2 - t)^2}{4z(z-1)(z-t)}$$
for some $|t|_v < 1$.  Observe that $\infty$ is a fixed point for $f$, and the points $a = 0, 1, t$ are sent to $\infty$ under one iteration of $f$.  In homogeneous coordinates, we have 
	$$F(z,w) = ((z^2 - t w^2)^2, 4zw(z-w)(z-tw)).$$

For simplicity, we will assume that $v(t) = 1$.  The general case has an identical proof.  Also, it is convenient to set $v(0) = \infty$ and $v(\infty) = -\infty$.  

\begin{prop}  \label{tent detail}
For any $r \in \Q$ and point $a\in \C_v$ with $v(a) = r$, we have 
	$$\sigma(F,a) = \left\{  \begin{array}{ll} 
			0 & \mbox{if } r \leq 0 \\
			2r & \mbox{if } 0 \leq r \leq 1 \\
			2 & \mbox{if } r \geq 1
			\end{array} \right.$$
and
	$$v(f(a)) \; \left\{  \begin{array}{ll} 
			< 0 & \mbox{if } r < 0 \\
			< 0 & \mbox{if } r = 0 \mbox{ and } r_v(a) = 1 \\
			=0 & \mbox{if } r = 0 \mbox{ and } r_v(a) \not=1 \\
			= 2r & \mbox{if } 0 < r < 1/2  \\
			= 1 & \mbox{if } r = 1/2 \mbox{ and } r_v\left(a/t^{1/2}\right) \not=\pm 1 \\
			> 1 & \mbox{if } r = 1/2 \mbox{ and } r_v\left(a/t^{1/2}\right) = \pm 1 \\
			= 2 - 2r & \mbox{if } 1/2 < r < 1  \\
			= 0 & \mbox{if } r = 1 \mbox{ and } r_v(a/t) \not= 1 \\
			< 0 & \mbox{if } r = 1 \mbox{ and } r_v(a/t) = 1 \\
			< 0 & \mbox{if } r > 1  
		\end{array}  \right.$$
where $r_v(\alpha)$ represents (as always) the residue class of $\alpha\in \O_v$.
\end{prop}

\begin{figure} [h]
\includegraphics[width=1.75in]{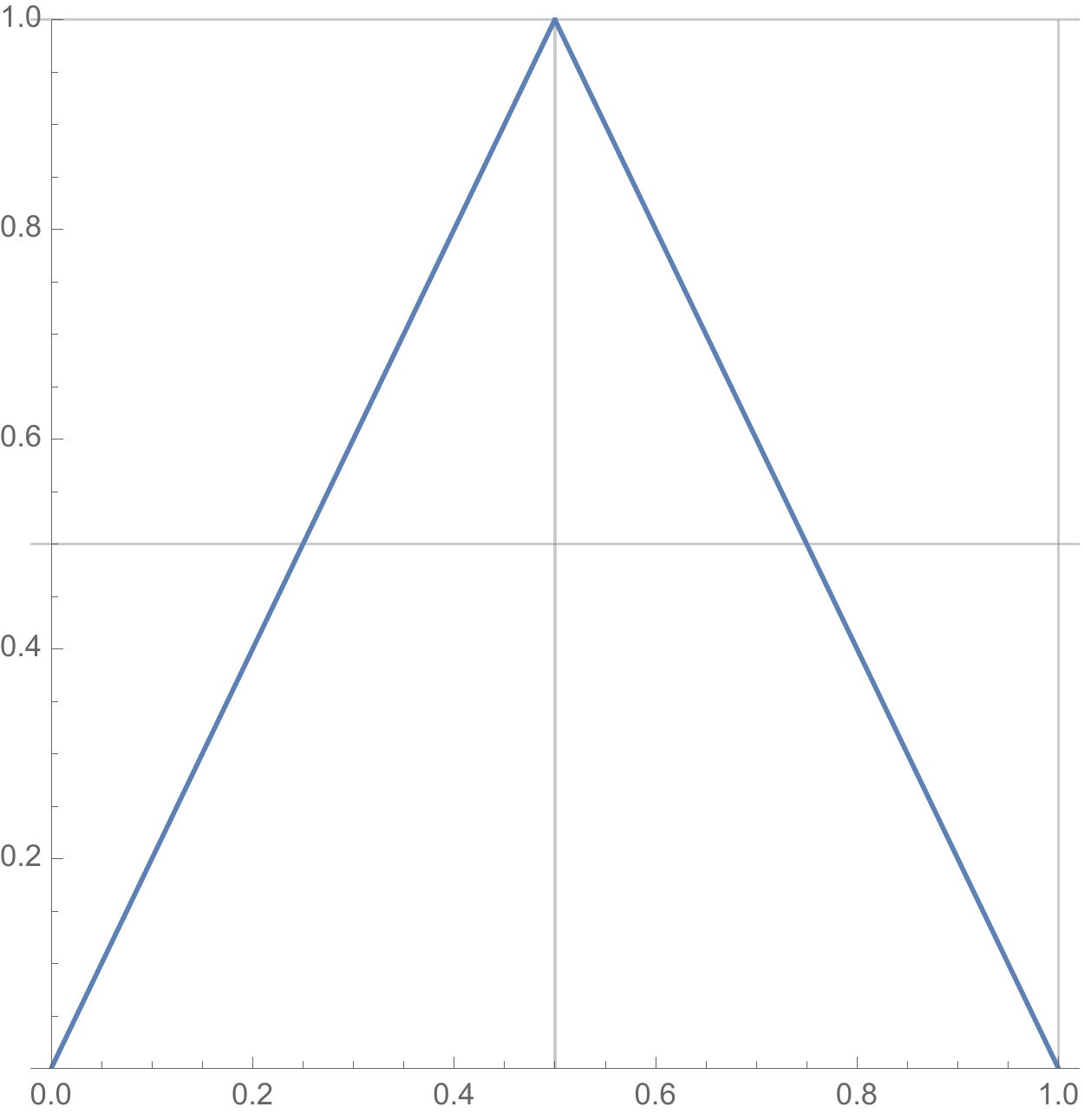}
\caption{ \small The graph of the tent map $T_2$ of slope $\pm 2$ on the interval $[0,1]$.}
\label{Tent}
\end{figure}

The first piece of Proposition \ref{tent detail} implies that the spine of $f$ is the closed interval 
	$$S_f = [\zeta_0, \zeta_1] \subset \bP^1_v,$$ 
where $\zeta_r$ is the Type II point representing the disk $D(0, |t^r|)$.  The second piece of the Proposition describes the action of $f$ on the spine.  It is exactly the tent map $T_2$ of slope $\pm 2$ on the unit interval; see Figure \ref{Tent}.  
This is also included in the statement of \cite[Proposition 5.1]{FRL:ergodic}.  For this example, the Julia set of $f$ coincides with the spine, because it is totally invariant for $f$.  The two pieces of the Proposition together allow us to iterate points and compute the sequence $c_n = \sigma(F, f^n(a))$ for any starting point $a$.  The key observation is that all rational points in the interval $[0,1]$ are preperiodic for the tent map.

\proof[Proof of Proposition \ref{tent detail}]
Recall that we are assuming $v(t) = 1$.  Fix $a \not=0, 1, t$ with $|a|_v\leq 1$, and write $A = (a, 1)$.  Let $r = v(a)$.  Then 
	$$F(A) = (a^4 - 2 t a^2 + t^2, -4a^2 + 4a^3  + 4 t a - 4 t a^2).$$
For $0 < r < 1$, the minimum valuation occurs only for the term $-4a^2$, so $\sigma(F,a) = 2r$.  For $r = 0$, the first coordinate has valuation 0 while the second coordinate has valuation $\geq 0$, so again, $\sigma(F, a) = 2r = 0$.  For $r= 1$, the first coordinate has valuation 2, while the second coordinate has valuation $\geq 2$, and again, $\sigma(F, a) = 2 r = 2$.  For $r > 1$, we see that the minimal valuation occurs only for $t^2$, so we have $\sigma(F, a) = 2$.  For $r < 0$, we should consider a lift $A = (1, b)$ with $v(b) = -r > 0$.  Then 
	$$F(A) = (1 - 2 b^2 t + b^4 t^2, 4b - 4b^2 - 4tb^2 + 4tb^3).$$
and we see immediately that the minimal valuation occurs only for the term $1$, so $\sigma(F, a) = 0$.

Assuming $0 \leq r \leq 1$, upon dividing $F(A)$ by $t^{2r}$, we see that the first coordinate of $F(A)/t^{2r}$ has order $=\min\{2-2r, 2r, 1\}$, unless $r = v(a) = 0, 1/2, 1$ when there can be cancellation.  This shows that $v(f(a))$ has desired value for $0 < r < 1/2$ and $1/2 < r < 1$.  

To understand the case $r = v(a) = 1/2$, we can compute directly that  
	$$f(z t^{1/2})/t = \frac{-(z+1)^2(z-1)^2}{4 z^2} + O(t^{1/2}).$$
The leading term in this expansion describes the action of $f$ on the space of directions from the point $\zeta_{1/2}$ to the point $\zeta_1$.  It is of degree 4, sending both 0 and $\infty$ to $\infty$, and sending $\pm 1$ to 0.  This proves the desired statement for $r = 1/2$.  

For $r = 0$, we compute that 
	$$f(z) = \frac{z^2}{4(z-1)} + O(t),$$
so we know that the Gauss point $\zeta_0$ is fixed.  The action on directions to $\zeta_0$ fixes 0 and $\infty$, and sends 1 to $\infty$, but all other directions will have $v(f(a)) = v(a) = 1$.  This proves the desired statement for $r = 0$.  We have also proved the statement for $r < 0$, because the direction of $\infty$ from $\zeta_0$ is fixed (and the Berkovich disk is forward invariant).  

For $r=1$, we compute that 
	$$f(z t) = \frac{-1}{4z(z-1)} + O(t),$$
showing that $f$ sends $\zeta_1$ to the Gauss point $\zeta_0$ and describing the action on directions.  Indeed, the two directions 0 and 1 are send to $\infty$, while $\infty$ is the unique zero.  We conclude that $v(a) = 1 \implies v(f(a)) = 0$ for all $z =a/t$ satisfying $r_v(z) \not= 1$, and $v(f(a)) < 0$ for all $z = a/t$ satisfying $r_v(z) =1$.  The computation also shows that $v(a) > 1$ implies that $v(f(a)) < 0$, completing the proof of the proposition.  
\qed

\subsection{Proof of Theorem \ref{dynamical E}}
Fix $a\in \P^1(\C_v)$ and set $c_n = \sigma(F, f^n(a))$.  Recall from Proposition \ref{rationality test} that it suffices to show that the sequence $c_n$ of rational numbers is eventually periodic.  Look at the sequence of rational numbers $r_n = v(f^n(a))$.  If $r_n \in [0,1]$ for all $n$, then from Proposition \ref{tent detail}, the sequence $r_n$ is equal to a rational orbit of the tent map $T_2$, so it must be preperiodic.  Then, since $c_n = 2 r_n$, we conclude that $c_n$ also cycles eventually.  If there exists an $n_0$ so that $r_{n_0} < 0$, then Proposition \ref{tent detail} implies that $r_n <0$ for all $n \geq n_0$, and $c_n = 0$ for all $n\geq n_0$.  If there exists an $n_0$ so that $r_{n_0} > 1$, then $r_{n_0+1} < 0$ and $c_n =0$ for all $n> n_0$.  We conclude that $\eta(F,a) = \sum_n c_n/d^n$ is a rational number.  Therefore, from Lemma \ref{local height}, the value $\lhat_{f,v}(a)$ is also rational.
\qed

\bigskip
\section{Irrationality}
\label{irrationality}

In this section, we prove Theorems \ref{irrational} and \ref{all reals}.  

Actually, we prove a more technical statement, about the realization of arbitrary sequences of two numbers as the order itinerary of a point in $\P^1(\C_v)$.  Let $\mathbb{N}_0$ denote the non-negative integers.

\begin{prop} \label{irrational general}
Suppose $f$ fixes the Gauss point $\zeta_0$ in $\bP^1_v$. 
Let $F$ be a normalized homogeneous presentation of $f$.  Suppose $f(\zeta)=\zeta_0$ where $\zeta$ is a Type II point in the direction from $\zeta_0$ containing 0.  Assume there exist:  
\begin{itemize}
\item a backward orbit $\{b_n\}\subset \P^1(\C_v)$ of 0 (satisfying $f(b_1)=0$ and $f(b_{n+1}) = b_n$ for all $n\geq 1$) 
such that 
	$$\sigma(F,b_n) = 0$$ 
for all $n$;  
\item  a point $x\in \P^1(\C_v)$ such that $\sigma(F, f^n(x))=0$ for all $n\geq 0$;
\item  preimages of $x$, $0$, and $b_n$ in directions $U$ from $\zeta$ where $\sigma(F, \cdot)|U = \sigma_F(\zeta)$.  
\end{itemize}
Then $c := \sigma_F(\zeta) > \sigma_F(\zeta_0) =0$, and for any sequence $s = s_0s_1s_2\cdots \in \{0, c\}^{\mathbb{N}_0}$, there exists a point $a\in \P^1(\C_v)$ for which 
	$$\sigma(F, f^n(a)) = s_n$$
for all $n\geq 0$.
\end{prop}

The idea of the proposition is to exploit a standard dynamical construction.  Roughly, we will find disjoint regions $D_0$ and $D_1$ in $\P^1(\C_v)$ so that $f$ maps each region over the union $D_0 \cup D_1$.  Then, appealing to the completeness of the field $\C_v$, we conclude that for any itinerary $s = (s_n) \in \{0,1\}^\N$, there will be a point $x \in D_{s_0}$ so that $f^n(x) \in D_{s_n}$ for all $n\geq 1$.  The hypotheses of Proposition \ref{irrational general} allow us to apply Theorem \ref{spine shape} to find suitable regions $D_0$ and $D_1$ so that $\sigma(F,\cdot)$ is constant on each $D_i$, with $\sigma(F, D_0) \not= \sigma(F, D_1)$.

\subsection{Proof of Theorem \ref{irrational}}
Assume Proposition \ref{irrational general}.  Fix $d\geq 2$ and $t \in k$ with $v(t) = 1$.  Choose a rational function $g \in K(z)$ of degree $d-1$ for which $0$ has infinite orbit and $g(\infty) = \infty$.  For example, we could take the polynomial $g(z) = z^{d-1}+1$.  Consider 
	$$f(z) = g(z) \; \frac{z-t}{z+t} \in k(z).$$
In homogeneous coordinates, let us write 
	$$F(z,w) = ((z-tw)P_g(z,w), (z+tw)Q_g(z,w))$$
where $g(z) = P_g(z,1)/Q_g(z,1)$, with $F$ normalized as in \eqref{general F}.

\begin{lemma}  \label{irrational spine}
The spine of $f$ is the interval $[\zeta_0, \zeta_1]$ where $\zeta_0$ is the Gauss point and $\zeta_1 = D(0,|t|)$, and $\sigma(F, a)=r$ for all $v(a)=r$, with $0 \leq r \leq 1$.  
\end{lemma}

\proof
The spine must lie in the direction of 0 from $\zeta_0$, because 0 is the only common root of the reduction of $F$.  Observe that 
	$$F(z,w) = (zP_g(z,w) + O(t), zQ_g(z,w) + O(t)).$$
Thus, for all $a$ with $0 < v(a) < 1$, we can take $A  = (a, 1)$, and we find that $\sigma(F, a) = v(a)$.  On the other hand, since there are no other zeros and poles of $f$ besides $t$ and $-t$ in the direction of 0 from $\zeta_0$, we use \eqref{small a} to show that $\sigma(F, a) = 1$ for all $a$ with $v(a) \geq 1$.  
\qed

\medskip
The Gauss point $\zeta_0$ is fixed by $f$, and the induced map on directions is given by $g$.  As the degree of $g$ is $d-1$, there is a unique other preimage of $\zeta_0$, and this point is $\zeta_1 = D(0,|t|_v)$, mapping with degree 1; in fact 
	$$f(t z) = g(0) \frac{z-1}{z+1} + O(t),$$
describing the induced map on directions from $\zeta_1$ to $\zeta_0$.  Since the point $t$ maps to 0 by $f$, the direction from $\zeta_1$ containing $t$ will map over $\zeta_1$.  More generally, the expression for $f(tz)$ implies that the direction from $\zeta_1$ containing $z t$ maps surjectively, and with degree 1, to the direction from $\zeta_0$ containing $g(0) \frac{z-1}{z+1}$, for all $z\in K$.  Moreover, for each $b$ with $|b|_v = 1$, the direction $U_b$ from $\zeta_0$ containing $b$ is sent to the direction $U_{g(b)}$ from $\zeta_0$, surjectively, with multiplicity equal to that of $g$ at $b$.

We now apply Proposition \ref{irrational general}.  We let $\{b_n\}$  be any backward orbit of 0 for $f$; it will satisfy the hypothesis of Proposition \ref{irrational general} because the spine is contained in the direction of 0.  We let $x = \infty$.  Therefore, for any sequence $s \in \{0,1\}^{\mathbb{N}_0}$, there exists a point $a\in \P^1(\C_v)$ for which 
	$$\sigma(F, f^n(a)) = s_n.$$
Now choose any aperiodic sequence.  It follows from Proposition \ref{rationality test} and Lemma \ref{local height} that the corresponding points will have irrational local canonical height.  This completes the proof of Theorem \ref{irrational}.

\subsection{Proof of Theorem \ref{all reals}}
Let $k = \C(t)$ and $v = \ord_{t=0}$ so that $v(t) = 1$.  Let 
	$$f(z) = \frac{(z+1)(z-t)}{z+t}.$$
Note that this $f$ is an example of the type given in the beginning of the proof of Theorem \ref{irrational}, above, with $d=2$.  Therefore, from the proof of Theorem \ref{irrational}, for any sequence $s \in \{0,1\}^{\mathbb{N}_0}$, there exists a point $a \in \P^1(\C_v)$ for which 
	$$\sigma(F, f^n(a)) = s_n.$$
In the proof above, we chose the point $x$ (which satisfies $\sigma(F, f^n(x)) = 0$ for all $n\geq 0$) to be $\infty$; let us take instead $x=1$.  In this way, we continue to have $\sigma(F, f^n(x)) = 0$ for all $n\geq 0$, but in addition we have that $v(x) = 0$.

Now fix any real number $\alpha \in [-1, 0]$.  Then $-\alpha$ has a binary expansion of the form 
	$$-\alpha = \frac{1}{2}  \sum_{n=0}^\infty \frac{c_n}{2^n}$$
with $c_n \in \{0,1\}$ for all $n\geq 0$.  Taking this sequence $\{c_n\}$ for the itinerary, obtain a point $a \in \P^1(\C_v)$ satisfying $v(a) \geq 0$ for which 
	$$\eta(F,a) = \frac{1}{2}  \sum_{n=0}^\infty \frac{\sigma(F, f^n(a))}{2^n} = \frac{1}{2}  \sum_{n=0}^\infty \frac{c_n}{2^n} = -\alpha$$
and therefore
	$$\lhat_{f,v}(a) = \alpha$$
by Lemma \ref{local height}.  This completes the proof of Theorem \ref{all reals}.

\subsection{Proof of Proposition \ref{irrational general}}
Without loss of generality, we may assume that $\zeta = D(0,r)$ for some $0 < r < 1$.  From Proposition \ref{order function} and Theorem \ref{spine shape}, we know that $\sigma_F(\zeta) > 0$.  Furthermore, for $b$ in all but finitely many directions from $\zeta_0$, we have $\sigma(F, b) = 0$.   

Notice that the direction of $\infty$ from $\zeta$ is necessarily a bad direction for $\zeta$, since $\zeta_0$ is fixed.  There are at most finitely many other bad directions.  Let $z_1, \ldots, z_N$ be points in each of these bad directions from $\zeta$ (if they exist).  

We now fix notation for certain subsets of $\P^1(\C_v)$ to be used for the remainder of the proof.  For each $b \in \C_v$ with $|b|_v = 1$, we let 
	$$X_b = \{z \in \C_v:  |b - z|_v < 1\},$$
containing all points in the direction of $b$ from $\zeta_0$.   Let 
	$$X_0 = \{z \in \C_v: |z|_v < 1\}$$
and
	$$X_\infty = \{z \in \C_v :  |z|_v >1\}.$$
We set
	$$X_\zeta = D(0,r) = \{ |z|_v \leq r\}$$
and
	$$X_\zeta^{\mathit{good}} = X_\zeta \setminus \bigcup_{j = 1}^N \{z: |z_j - z|_v < r\}.$$
By hypothesis, we have $\sigma(F, X_{b_n}) = 0$ for all $n\geq 1$, and $\sigma(F, X_{f^n(x)}) = 0$ for all $n\geq 0$.  It follows from Theorem \ref{spine shape} that none of the $b_n$ or $f^n(x)$ lie in bad directions from $\zeta_0$, and therefore $f(X_{b_{n+1}}) = X_{b_n}$ and $f(X_{f^{n-1}(x)}) = X_{f^n(x)}$ for all $n\geq 1$ and $f(X_{b_1}) = X_0$.  We also know that $\sigma(F, X_\zeta^\mathit{good}) = \sigma_F(\zeta)$.  Furthermore, there exist preimages $b_n'$ of $b_n$ and $x'$ of $x$ in $X_\zeta^{\mathit{good}}$; we set
	$$Y_n = \{|z - b_n'|_v < r\}$$
and 
	$$Y_x = \{|z - x'|< r\}$$
so that $\sigma(F, Y_n) = \sigma(F, Y_x) = \sigma_F(\zeta)$ for all $n$.  

We now begin our construction.  

Set $c = \sigma_F(\zeta)$.  Fix sequence $s \in \{0,c\}^{\mathbb{N}_0}$.  We will construct a nested sequence of regions $D_n\subset \P^1(\C_v)$, for $n\geq 0$, for which 
	$$\sigma(F, f^j(a)) = s_j$$
for all $j\leq n$ and all $a \in D_n$.  The conclusion is obtained by choosing any point 
	$$a \in \bigcap_{n\geq 0} D_n$$
which we will show is nonempty, appealing to the completeness of $\C_v$.   

We proceed by induction.  For the base case, if $s_0 = c$, we set $D_0 = X_\zeta^{\mathit{good}}$.  Then 
	$$\sigma(F,a) = c = s_0$$ 
for all $a\in D_0$.  If $s_0 = 0$, we let 
	$$m(0) = \sup\{n \geq 0: s_0 = s_1 = \cdots = s_n = 0\}.$$
For $m(0) = \infty$, we take $D_j = X_x$ for all $j \geq 0$.  Then 
	$$\sigma(F, f^j(a)) = 0 = s_j, \; \forall j \geq 0, \forall a \in X_x$$
and we are done.  On the other hand, if $m(0)$ is finite, we set $D_0 = D_1 = \cdots = D_{m(0)} = X_{b_{m(0)+1}}$.  Note that $f^{m(0)+1}(X_{b_{m(0)+1}}) = X_0 \supset X_\zeta$.  Thus, there exists a (closed) disk $U_{m(0)+1} \subset D_{m(0)}$ of radius $\leq r$ sent by $f^{m(0)+1}$ to $X_\zeta$, and a region $D_{m(0)+1} \subset U_{m(0)+1} \subset D_{m(0)}$ sent onto $X_\zeta^{\mathit{good}}$ by $f^{m(0)+1}$.  By construction, we have 
	$$\sigma(F,a) = \cdots = \sigma(F, f^{m(0)}(a)) = 0 \mbox{ and } \sigma(F, f^{m(0)+1}(a)) = c$$
for all $a\in D_{m(0)+1}$, so that 
	$$\sigma(F, f^j(a)) = s_j, \; \forall j \leq m(0)+1, \; \forall a \in D_{m(0)+1},$$
and $f^{m(0)+1}(D_{m(0)+1}) = X_\zeta^{\mathit{good}}$.  

We now explain the induction step.  Assume that we have a region $D_k$ so that 
	$$\sigma(F,f^j(a)) = s_j$$
for all $j\leq k$ and all $a\in D_k$, with $s_k = c$ and $f^k(D_k) = X_\zeta^{\mathit{good}}$.  Assume further that the region $D_k$ is either an open disk of radius $r_k>0$ or it is the complement of finitely many open disks in a closed disk of radius $r_k$.  Suppose $s_{k+1} = c$.  Let $D$ be a disk inside $X_\zeta^{\mathit{good}}$ sent by $f$ to $X_\zeta$.  Then we let $U_{k+1}$ be the disk inside $D_k$ sent by $f^k$ to $D$; finally, we let $D_{k+1}$ be the subset of $U_{k+1}$ sent by $f^{k+1}$ to $X_\zeta^{\mathit{good}}$.  Note that  $D_{k+1}$ is itself the complement of finitely many open disks in $U_{k+1}$, which is a disk of radius $\leq r r_k$.  Moreover, $\sigma(F, f^{k+1}(a)) = c = s_{k+1}$ for all $a\in D_{k+1}$, with $f^{k+1}(D_{k+1}) = X_\zeta^{\mathit{good}}$, completing the induction step in this case. 

On the other hand, for $s_{k+1} = 0$, we let $m(k+1)\in \N\cup\{\infty\}$ be defined as 
	$$m(k+1) := \sup\{n \geq k+1: s_{k+1} = s_{k+2} = \cdots = s_n = 0\}$$
For $m(k+1)=\infty$, we set $D = Y_x$, sent to $X_x$ by $f$, and let $D_{k+1}$ be an open disk inside $D_k$ which is sent by $f^k$ to $D$.  We put $D_j = D_{k+1}$ for all $j\geq k+1$.  Then 
	$$\sigma(F, f^n(a)) = s_n \; \forall n$$
and we are done.  But if $m(k+1)$ is finite, we set $D = Y_{m(k+1)-k}$, which is sent by $f$ to $X_{b_{m(k+1)-k}}$.  Then we let $D_{k+1} = D_{k+2} = \cdots = D_{m(k+1)}$ be the preimage of $D$ by $f^k$ inside $D_k$, so that $D_{k+1}$ is a disk of radius $r_k$ with $f^{k+1}(D_{k+1}) = X_{b_{m(k+1)-k}}$ and $f^{m(k+1)+1}(D_{k+1}) = X_0 \supset X_\zeta$.  Then finally, we let $U_{m(k+1)+1}\subset D_{m(k+1)}$ be the preimage of $X_\zeta$, so $U_{m(k+1)+1}$ is a disk of radius $\leq r r_k$, and we let $D_{m(k+1)+1}$ be the region inside $D_{k+1}$ sent by $f^{m(k+1)+1}$ to $X_\zeta^{\mathit{good}}$.  Then $D_{m(k+1)+1}$ is the complement of finitely many open disks in $U_{m(k+1)+1}$, with
	$$\sigma(F, f^j(a)) = 0 = s_j \; \forall \; k+1 \leq j \leq m(k+1),$$ 
and
	$$\sigma(F, f^{m(k+1)+1}(a)) = c = s_{m(k+1)+1}$$
for all $a\in D_{m(k+1)+1}$, and $f^{m(k+1)+1}(D_{m(k+1)+1}) = X_\zeta^{\mathit{good}}$.  This completes the induction argument.  

It remains to observe that the nested intersection $\bigcap_n D_n$ is nonempty.  If the series $s$ ends with an infinite string of 0s, then the nested intersection coincides with one of the entries $D_n$, which is a nonempty open disk by construction.  In general, when there are infinitely many $n$ for which $s_n = c$, we appeal to the completeness of the field $\C_v$.  By the induction argument, for each $n \geq 0$, there exists a disk $U_n$ such that $D_n \supset U_{n+1} \supset D_{n+1}$ and the radii of the disks $U_n$ shrink by at least a factor of $r$ for each $n$ with $s_n = c$.  Then 
	$$\bigcap_n D_n = \bigcap_n U_n$$
is nonempty by the completeness of $\C_v$.  

This completes the proof of Proposition \ref{irrational general}.

\bigskip
\section{Quadratic maps}
\label{degree 2}

In this section, we prove Theorem \ref{quadratic}.  We begin by stating a result of Kiwi in \cite{Kiwi:quad}.  
Applying Kiwi's theorem, our proof will break into four cases.  

\subsection{Julia sets in degree 2}
Suppose $f$ is a quadratic rational function defined over $k = K(X)$ (where $X$ is a smooth curve defined over $K$).  Fix $v \in M_k$.  Kiwi classified the action of $f$ on its Julia set in the Berkovich line $\bP^1_v$.  We do not use the full strength of his results.  Here is a reformulation of his theorems that will suit our particular needs, combining the main results of \cite{Kiwi:quad} with some preliminary results appearing in \cite{Milnor:quad, Epstein:bounded, D:moduli2}. 

\begin{theorem}  \label{Kiwi}
Suppose $f$ has degree 2 and does not have potential good reduction at $v$.  Then exactly one of the following holds:
\begin{enumerate}
\item  $f$ has an attracting fixed point.  In this case, the Julia set $J(f)$ is a Cantor set in $\P^1(\C_v)$; all points of $\bP^1_v \setminus J(f)$ tend to the attracting fixed point under iteration.
\item  $f$ has a fixed Type II point on which, in suitable coordinates, $f$ acts on the space of directions by 
	$$z \mapsto \lambda z$$
for $\lambda \in K^*$ not a root of unity, or by
	$$z \mapsto z + 1,$$
and the unique bad direction points towards $z=1$. 
\item $f$ has a fixed Type II point $\zeta$ on which, in suitable coordinates, $f$ acts on directions by 
	$$z \mapsto \lambda z$$
for $\lambda \in K^*$ a primitive $p$-th root of unity, for some $p>1$; and in the direction of $z=1$ from $\zeta$, $f$ has a periodic Type II point of period $p$ on which, in suitable coordinates, $f^p$ acts on directions by 
	$$w \mapsto w + \tau + 1/w$$
for some $\tau \in K$, with exactly two bad directions pointing towards $w=0, \infty$ and with $\zeta$ lying in the direction of $w=\infty$. 
\item  $f$ has a fixed Type II point $\zeta$ on which, in suitable coordinates, $f$ acts on directions by 
	$$z \mapsto \lambda z$$
for $\lambda \in K^*$ a primitive $p$-th root of unity, for some $p>1$; and in the direction of $z=1$ from $\zeta$, $f$ has a periodic Type II point of period $p$ on which, in suitable coordinates, $f^p$ acts by 
	$$w \mapsto w + 1$$
with exactly two bad directions pointing towards $w = 0, \infty$ and with $\zeta$ lying in the direction of $w=\infty$.  
\end{enumerate}
\end{theorem}

\noindent
The history of this statement begins with Milnor's description of the moduli space $\M_2$ of quadratic maps over $\C$, in \cite[\S3]{Milnor:quad}.  He observed that $\M_2 \iso \C^2$ via the identification of the conjugacy class of $f$ with (two of) the symmetric functions in the three fixed-point multipliers.  Milnor defined a compactification $\Mbar_2 \iso \P^2$ where the boundary points are identified with (unordered) triples $\{\lambda, 1/\lambda, \infty\}$ of limiting fixed-point multipliers, with $\lambda \in \C$; for $\lambda\not=0$, these boundary points correspond to the action of $f$ on directions at the fixed Type II points appearing in cases (2), (3), and (4).  Epstein studied the reduction maps of degree 2 exhibited in case (3) to analyze the structure of hyperbolic components in $\M_2$ \cite{Epstein:bounded}.  Such nontrivial reduction maps at periodic Type II points are called ``rescaling limits" and arise as boundary points of $\M_2$ in more elaborate compactifications; see \cite{Kiwi:rescaling} and \cite{D:moduli2}.  

\begin{remark}
The example of Theorem \ref{all reals} satisfies condition (2) of Theorem \ref{Kiwi}, with the reduction map $z\mapsto z+1$.
\end{remark}

\subsection{Fixed-point multipliers}  \label{multiplier}
The {\em multiplier} of a fixed point of $f$ is simply the derivative of $f$ at that point.  The four cases of Theorem \ref{Kiwi}, and the value of $\tau$ in case (3), can be characterized in terms of the three fixed-point multipliers $\{\alpha, \beta, \gamma\}$ of a quadratic map.  To see this, first recall that the conjugacy class of a quadratic map $f$ is uniquely determined by its triple of fixed-point multipliers.  In addition, the multipliers are subject to exactly one condition, that 
\begin{equation} \label{fixed point formula}
	\alpha \beta \gamma - (\alpha + \beta + \gamma) + 2 = 0
\end{equation}
\cite[Lemma 3.1]{Milnor:quad}; see also \cite[Theorems 1.14 and 4.56]{Silverman:dynamics}.  For a characterization of potential good reduction in degree 2, we have:

\begin{prop} \label{quadratic reduction}
Fix $f$ of degree 2 defined over $k$, and fix a place $v$ of $k$.  Then $f$ has potential good reduction at $v$ if and only if the three fixed-point multipliers $\{\alpha, \beta, \gamma\}$ satisfy
	$$|\alpha|_v, |\beta|_v, |\gamma|_v \leq 1.$$
\end{prop}

\proof
Suppose that $f$ has potential good reduction.  Then under a M\"obius change of coordinates, $f$ will have good reduction.  For any map of good reduction, the multipliers at all periodic points will be $\leq 1$ in absolute value \cite[Chapter 2]{Silverman:dynamics}.  As the multipliers are invariant under coordinate change, we deduce that $|\alpha|_v, |\beta|_v, |\gamma|_v \leq 1$.  Conversely, suppose the three fixed point multipliers satisfy $|\alpha|_v, |\beta|_v, |\gamma|_v \leq 1$.  Recalling that $r_v$ denotes the reduction map modulo $v$, assume that $r_v(\alpha) \not=1$, so that $r_v(\alpha\beta) \not=1$ by \eqref{fixed point formula}.  Then, exactly as in the proof of \cite[Lemma 3.1]{Milnor:quad}, we observe that in suitable coordinates on $\P^1(\C_v)$ the map $f$ can be expressed explicitly in terms of its multipliers.  For example, we may write 
	$$f(z) = \frac{z^2 + \alpha z}{\beta z + 1},$$
which has a fixed point of multiplier $\alpha$ at $z=0$ and a fixed point of multiplier $\beta$ at $\infty$.
The resultant of a normalized homogeneous presentation $F$ of $f$ is computed as $\Res(F) = 1 - \alpha\beta$ so that $|\Res(F)|_v = 1$.   That is, in these coordinates, $f$ will have good reduction.   See also \cite[p. 193]{Silverman:dynamics}.  If, on the other hand, we have $r_v(\alpha) = r_v(\beta) = r_v(\gamma) = 1$, then we may choose coordinates on $\P^1(\C_v)$ so that the two critical points of $f$ lie at 1 and $-1$ and the point $0$ is fixed with multiplier $\alpha$.  Then we have
	$$f(z) = \frac{\alpha z}{z^2 + \delta z + 1}$$
for some $\delta \in \kbar$, which will have good reduction if and only if $|\delta|_v \leq 1$.  To determine $\delta$ in terms of $\{\alpha, \beta,\gamma\}$, we first compute that the other two fixed-point multipliers of $f$ are (up to ordering)
	$$\beta = (4 - 2 \alpha - \delta^2 + \delta \sqrt{-4 + 4 \alpha + \delta^2})/(2 \alpha)$$
and
	$$\gamma = (4 - 2 \alpha - \delta^2 - \delta \sqrt{-4 + 4 \alpha + \delta^2})/(2 \alpha).$$
Eliminating $\alpha$ from these equations shows that
	$$\beta^2 + \gamma^2 + \beta \gamma (\delta^2 - 2) - \delta^2 = 0.$$
Solving for $\delta$, we find
	$$\delta = \pm \frac{\beta - \gamma}{\sqrt{1 - \beta\gamma}}$$
as long as $\beta\gamma \not= 1$.  For $\beta \gamma = 1$, so that $\beta = \gamma = 1$ by \eqref{fixed point formula}, we compute directly from the formulas for $\beta$ and $\gamma$ that 
	$$\delta = \pm 2 \sqrt{1 - \alpha}.$$
In either case, recalling that $r_v(\alpha) = r_v(\beta) = r_v(\gamma) = 1$, we have $|\delta|_v < 1$, and we conclude that $f$ has good reduction at $v$.  
\qed

\medskip
It follows from Proposition \ref{quadratic reduction} that, in the absence of potential good reduction, there will be at least one multiplier with absolute value $>1$.  The first case of Theorem \ref{Kiwi} is already formulated in terms of the fixed points.  If there is an attracting fixed point, then the other two multipliers must have absolute value $>1$, by \eqref{fixed point formula}.   The remaining cases of Theorem \ref{Kiwi} hold when two of the multipliers, say $\alpha$ and $\beta$, have absolute value 1, and their residue classes are given by $r_v(\alpha) = \lambda$ and $r_v(\beta) = 1/\lambda$.  The value of $\lambda$ will determine whether we are in case (2) or in cases (3) and (4).  Cases (3) and (4) are distinguished by the residue class
\begin{equation} \label{tau squared}
  \tau^2 := r_v\left(\frac{\alpha^p - 1}{1 - \alpha\beta}\right) = r_v\left(\frac{\beta^p - 1}{1 - \alpha\beta}\right).
\end{equation}
For $\tau^2\ne \infty$, then the value of $\tau$ in case (3) is either choice of square root of $\tau^2$; otherwise we are in case (4).  See \cite[\S5]{D:moduli2} and \cite[Proposition 6.1, Proposition 6.4]{D:moduli2}.

\subsection{Attracting fixed point}
Now we work towards the proof of Theorem \ref{quadratic}, which will be completed in \S\ref{quadratic proof}.  Suppose that $f$ has an attracting fixed point, and assume that it does not have potential good reduction.  From Theorem \ref{Kiwi}(1), the Julia set of $f$ is a Cantor set in $\P^1(\C_v)$.  From Proposition \ref{quadratic reduction} and equation \eqref{fixed point formula}, the other two fixed points have multipliers of absolute value $>1$.  Suppose the three fixed point multipliers are $\alpha, \beta, \gamma$ with $|\gamma|_v < 1$ and $|\alpha|_v \geq |\beta|_v > 1$.  Again from \eqref{fixed point formula}, we compute that $|\alpha|_v > |\beta|_v = 1/|\gamma|_v >1$ if $|\alpha\gamma|_v > 1$; otherwise, we have $|\beta|_v = |\alpha|_v$.  We may assume we have chosen coordinates so the attracting fixed point lies at $\infty$.  

In  the proof of Theorem \ref{quadratic}, we use the following result:

\begin{prop}  \label{attracting fixed point}
Suppose $f$ has an attracting fixed point at $\infty$ and two repelling fixed points in $\P^1(\C_v)$ with multipliers $\alpha$ and $\beta$.  Then the local canonical height $\lhat_{f,v}$ satisfies 
	$$\lhat_{f,v}(a) \in \Q \quad \forall \; a\in \P^1(\C_v)\setminus\{\infty\}$$
if and only if $|\alpha|_v = |\beta|_v$.  
\end{prop}

For convenience, we will assume that $0$ is a repelling fixed point.  Assume further that the other zero of $f$ is at $1$.  We can express our map in homogeneous coordinates as 
	$$F(z,w) = (c \, z (z - w), \, u \, w (p_2 z - p_1 w))$$
where $\min\{v(c), v(u)\} = \min\{v(p_1), v(p_2)\} = 0$, with $p = [p_1:p_2]$ being the other pole of $f$ in $\P^1(\C_v)$.  We compute that the attracting fixed point at $\infty$ has multiplier $u p_2/c$, and therefore 
	$$|p_2u|_v < |c|_v.$$
Since 0 is repelling with multiplier $c/(p_1u)$, we deduce that
	$$|p_1u|_v < |c|_v.$$
Independent of $p$, then, we see that $|u|_v < |c|_v$, so we may take $c=1$ and $|u|_v < 1$.  Thus, we may write
\begin{equation} \label{special F}
	F(z,w) = (z (z - w), \, u \, w (p_2 z - p_1 w))
\end{equation} 
where $v(u)>0$ while $\min\{v(p_1), v(p_2)\} = 0$.  Computing directly, the triple of fixed point multipliers at the three fixed points $\infty$, 0, and $(1 - u p_1)/(1 - u p_2) \approx 1$ is 
$$\left\{u p_2, \; \frac{1}{u p_1},  \;\frac{1 - 2 u p_1 + u^2 p_1 p_2}{u(p_2 - p_1)} \right\}$$

\begin{lemma}  \label{pole location}
The two repelling fixed-point multipliers have the same absolute value if and only if $v(p) = v(p-1)$, where $p = p_1/p_2$.  
\end{lemma}

\proof
For $v(p)<0$, then necessarily $v(p-1) = v(p)$, and we can take $p_1 = 1$ and $|p_2|_v< 1$.  In this case both repelling multipliers have absolute value $= 1/|u|_v$.  For $v(p)= 0$, then we can again take $p_1 = 1$, and the condition that the two multipliers have the same absolute value $= 1/|u|_v$ is equivalent to $|p_2 - p_1|_v = 1$, which is equivalent to $v(p-1)=0$.  For $v(p) >0$, so that $v(p-1) = 0 \not= v(p)$, taking $p_2 = 1$ we see from the formulas that the repelling multipliers will have absolute values $1/|u p_1|_v > 1/|u|_v$.  
\qed

\medskip	
Now take any point $a \in \P^1(\C_v)$.  We apply Theorem \ref{spine shape} and Proposition \ref{order function} to the $F$ of \eqref{special F} to deduce that $\sigma(F,a) = 0$ unless $v(a)>0$ or $v(a-1)>0$.  From (\ref{small a}), we compute that 
\begin{equation} \label{polylike sigma}
	\sigma(F, a) = \min\{v(a) + v(a - 1), \; v(u) + \max\{0, v(a-p)\} \}
\end{equation}
for all $|a|_v \leq 1$, where $p = p_1/p_2$. The term involving $p$ will vanish if $|p|_v > 1$.   The shape of the spine will depend on the location of $p$ and the value of $u$.  Nevertheless:

\begin{lemma} \label{Julia disks}
There exist open disks $D_0 \ni 0$ and $D_1 \ni 1$ such that $f(D_i)\supset D_0\cup D_1$ and $\sigma(F, \cdot)$ is constant on $D_i$ for $i = 0,1$.  A conjugacy from $f|J(f)$ to the shift map on $\{0,1\}^{\mathbb{N}_0}$ is given by the itinerary of points in these two disks.  The two repelling fixed-point multipliers have the same absolute value if and only if $\sigma(F, D_0) = \sigma(F, D_1)$.  
\end{lemma}

\proof
First suppose that the two fixed point multipliers have the same absolute value.  From Lemma \ref{pole location}, we know that $v(p)= v(p-1)$.  Assume that $v(p) = 0$, and let $r = |u|_v < 1$. Then we can compute directly that the two preimages of the Gauss point $\zeta_0$ are $\zeta_r := D(0, r)$ and $\zeta_{1,r} := D(1, r)$.  Indeed, taking $p_1 = p$ and $p_2 = 1$, we have
	$$f(z) = \frac{z(z-1)}{u (z-p)},$$
so that
	$$f(z \, u) = z/p + O(u)$$
and
	$$f(1 + z \, u) = z/(1-p) + O(u).$$
It follows that the Julia set is the nested intersection of the preimages of the disk $D(0,1)$; it lies in the union of the open disks $D_0 = \{|z|_v < r\}$ and $D_1 = \{|z-1|_v < r\}$; and the itinerary of an orbit with respect to $D_0$ and $D_1$ determines the conjugacy from $f|J(f)$ to the 2-shift.  It is straightforward to see from (\ref{polylike sigma}) that 
	$$\sigma(F, a) = r = |u|_v$$
for all $a \in D_0 \cup D_1$.  

The case where $v(p) = v(p-1) < 0$ is similar.  In this case, we set $p_1 = 1$ and $p_2 = 1/p$ so that 
	$$f(z) = \frac{z(z-1)}{u(z/p - 1)} = \frac{pz (z-1)}{u (z-p)}.$$
We compute that 
	$$f(z \, u) = z + O(u)$$
and
	$$f(1 + z \, u) = p z/(1-p) + O(u),$$
so that, again, the preimages of the Gauss point $\zeta_0$ are the Type II points $\zeta_r$ and $\zeta_{1,r}$.  As before, we have 
	$$\sigma(F, a) = r = |u|_v$$
for all $a \in D_0 \cup D_1$.

Now assume the two fixed point multipliers do not have the same absolute value.  From Lemma \ref{pole location}, either $v(p)>0$ or $v(p-1)>0$.  Suppose that $v(p) >0$.  From (\ref{polylike sigma}), we compute that 
	$$\sigma(F, a) = v(u)$$
for all $a \in \{z: |z-1|_v < |u|_v\}$.  On the other hand
	$$\sigma(F, a) = v(u) + v(p) > v(u)$$
for all $a \in \{z: |z|_v < |u p|_v\}$.  Moreover, the closed disk $D(0, |u p|_v)$ is sent to the Gauss point $\zeta_0$, since 
	$$f(up z) = z + O(u)$$
and the disk $D(1, |u|_v)$ is also sent to $\zeta_0$, since
	$$f(u z) = z + O(p).$$
It follows that the Julia set lies in the union of $D_1 = \{z: |z-1|_v < |u|_v\}$ and $D_0 = \{z: |z|_v < |u p|_v\}$.  The Julia set lies in the union of the two disks $D_0$ and $D_1$, with the orbit itineraries determined by these disks providing the conjugacy to the 2-shift.   However, now we have that $\sigma(F, D_0) = v(u) + v(p) > v(u) = \sigma(F, D_1)$.  The case where $v(p-1)>0$ is similar.  
\qed

\medskip
\proof[Proof of Proposition \ref{attracting fixed point}]
For $a \not\in J(f)$, we have $f^n(a) \to \infty$ under iteration by Theorem \ref{Kiwi}(1). Thus, $f^n(a)$ will lie in the direction of $\infty$ from the Gauss point $\zeta_0$ for all $n>>0$, and therefore, from Theorem \ref{spine shape} we have $\sigma(F, f^n(a)) = 0$ for all $n >> 0$.  Thus, the sum
	$$\eta(F, a) = \frac{1}{2} \sum_{n=0}^\infty \frac{\sigma(F, f^n(a))}{2^n}$$
is rational, and $\lhat_{f,v}(a)$ will be rational by Lemma \ref{local height}.  

For $a \in J(f)$, we have $f^n(a) \in J(f)$ for all $n$.   From Lemma \ref{Julia disks}, we know that $\sigma(F, \cdot)$ is constant on $J(f)$ if the two repelling fixed-point multipliers have the same absolute value.  In that case, $\eta(F,a)$ is clearly rational, and Lemma \ref{local height} shows that $\lhat_{f,v}(a)$ will be rational.   On the other hand, also from Lemma \ref{Julia disks}, if $\sigma(F, \cdot)$ is nonconstant, then it takes distinct rational values on the two disks $D_0$ and $D_1$ that determine the conjugacy to the 2-shift.  Let $a$ be any point in $J(f)$ with an aperiodic itinerary with respect to the disks $D_0$ and $D_1$ of Lemma \ref{Julia disks}; then the sequence $\sigma_n = \sigma(F, f^n(a))$ is not eventually periodic.  Consequently, the value of $\eta(F,a)$ will be irrational.  (We may apply our rationality criterion of Proposition \ref{rationality test} after transforming the values by $\sigma_n \mapsto (\sigma_n - \sigma(F, D_0))/(\sigma(F,D_1)- \sigma(F, D_0))$ so the terms in the series are all 0 and 1.)  We conclude from Lemma \ref{local height} that  $\lhat_{f,v}(a)$ is irrational.  
\qed

\subsection{Proof of Theorem \ref{quadratic}.}  \label{quadratic proof}
Fix a quadratic map $f$ and a place $v$ of $k$, and fix a normalized homogeneous presentation $F$ of $f$.  

Suppose first that $f$ has potential good reduction at $v$.  It follows from Proposition \ref{potential good reduction rationality} that the local canonical height of $f$ takes only rational values.

Assume now that $f$ does not have potential good reduction at $v$.  Then we are in one of the four cases of Theorem \ref{Kiwi}.  Case (1) of Theorem \ref{Kiwi} is treated by Proposition \ref{attracting fixed point}, showing that if there exists an attracting fixed point, then all points have rational local canonical height if and only if the two repelling fixed points have multipliers of the same absolute value.  In other words, all points have rational canonical height if and only if $f$ is strongly polynomial-like.

It remains to show that if we are in cases (2), (3) or (4) of Theorem \ref{Kiwi}, then there must exist at least one point $a\in \P^1(\C_v)$ for which $\lhat_{f,v}(a)$ is an irrational number.  

Assume that we are in case (2) of Theorem \ref{Kiwi}.  The result will follow from Proposition \ref{irrational general}.   Choose coordinates so that the fixed Type II point of Theorem \ref{Kiwi}(2) is the Gauss point $\zeta_0$.  From Theorem \ref{Kiwi}(2), there is a unique preimage $\zeta \not= \zeta_0$ of the Gauss point $\zeta_0$ in the direction of $z=1$ from $\zeta_0$.  The action on directions from $\zeta_0$ sends 1 to $z = \lambda$ (respectively, to $z=2$ when the reduction map is $z+1$); it follows that the direction from $\zeta$ containing $\zeta_0$ must also be sent towards $z=\lambda$ (respectively, towards $z=2$); otherwise, by continuity, there would be another preimage of $\zeta_0$ lying between the two points, which is impossible.  Thus, the map satisfies all the hypotheses of the Proposition \ref{irrational general}, taking the point 1 to play the role of 0, taking $b_n$ in the direction of $\lambda^{-n}$ (respectively, $1-n$) for all $n>0$, and setting $x = \infty$.  Theorem \ref{spine shape} shows that the order function $\sigma(F, \cdot)$ has the needed properties for Proposition \ref{irrational general}.  Therefore, choosing any sequence $s \in \{0, \sigma_F(\zeta)\}^\N$ which is not eventually periodic, we obtain a point $a\in \P^1(\C_v)$ with the given order itinerary $\sigma(F, f^n(a)) = s_n$.  It follows that 
	$$\eta(F,a) = \frac{1}{2} \sum_{n=0}^\infty \frac{\sigma(F, f^n(a))}{2^n}$$ 
will be irrational by Proposition \ref{rationality test}, and therefore also the local height by Lemma \ref{local height}.

If $f$ falls into cases (3) or (4), choose coordinates so that the periodic Type II point of period $p>1$ in Theorem \ref{Kiwi} is the Gauss point $\zeta_0$ with the coordinate $w$ on the space of directions shown in the theorem.  We work with the iterate $f^p$ and the distinguished bad direction $w=0$ from $\zeta_0$.  Unfortunately, we cannot apply Proposition \ref{irrational general} directly, because we may not be able to show that all the hypotheses are satisfied.  Nevertheless, the arguments in the proof of Proposition \ref{irrational general} may be repeated to show there exists at least one point with an aperiodic itinerary.  Specifically, observe that the spine of $f^p$ lies in the directions of 0 and $\infty$ from $\zeta_0$, because these are the unique bad directions by Theorem \ref{Kiwi}.  Take any end $\zeta \not= \zeta_0$ of the spine for $f^p$ lying in the direction of 0 from $\zeta_0$, so that $f^p(\zeta) = \zeta_0$.  (The existence of $\zeta$ is guaranteed by Theorem \ref{spine shape}.)  

Let $c = \sigma_{F^p}(\zeta) > 0$.  

Importantly, since $\zeta$ is an end of the spine, the order function $\sigma_{F^p}$ will be constant (equal to $c$) in all but one direction from $\zeta$.  Let $U$ denote this direction from $\zeta$ containing $\zeta_0$.  Consider the map of directions from $\zeta$ to $\zeta_0$; let $V$ denote the direction from $\zeta_0$ corresponding to the image of $U$.  Choose an infinite backward orbit $b_n$ of 0 satisfying the conditions of Proposition \ref{irrational general} for $f^p$.  The forms of the reduction map, $\bar{f^p}(w) = w + \tau + 1/w$ or $\bar{f^p}(w) = w+1$ presented in Theorem \ref{Kiwi} (3) and (4), show that such infinite orbits exist.  Take any point $x$ with forward orbit under $\bar{f^p}$ that does not contain 0 or $\infty$; in fact, we can always take either $x = 1$ or $x=-1$.  This $x$ will satisfy the hypotheses of Proposition \ref{irrational general}.  

Note that the direction $V$ contains at most one element of $\{0, x, b_1, b_2, b_3, \ldots\}$.  If $V$ contains none of these points, then we can apply Proposition \ref{irrational general} directly to conclude that for any sequence $s = (s_n) \in \{0,c\}^\N$, there exists a point $a$ with order itinerary 
	$$\sigma(F^p, f^{np}(a)) = s_n.$$
Suppose $V$ contains $x$.  Let $s = (s_n) \in \{0,c\}^\N$ be any aperiodic sequence.  Such sequences cannot end with an infinite string of 0s, so the proof of Proposition \ref{irrational general} goes through to show that for any such sequence, there exists a point $a$ with the given order itinerary for $f^p$.  Suppose instead that $V$ contains 0.  The proof of Proposition \ref{irrational general} shows that we can realize any sequence in $\{0,c\}^\N$ with no repeated $c$'s.  Finally, suppose that $V$ contains $b_n$.  The construction shows that we can realize any sequence in $\{0, c\}^\N$ that does not have strings of exactly $n$ zeroes that are preceded and followed by a $c$.  Thus, in every case, there exist points with aperiodic order itineraries and therefore points with irrational local canonical height for $f^p$, applying Proposition \ref{rationality test} and Lemma \ref{local height} to $f^p$. 

It remains to observe that a point with irrational local canonical height for an iterate $f^p$ will also be of irrational local height for $f$.  But this follows  from the computation of \eqref{order iterate}, where we observed that
	$$\sigma(F^p, a) = d^{p-1}\, \sigma(F, a) + d^{p-2} \, \sigma(F, f(a)) + \cdots + \sigma(F, f^{p-1}(a)).$$
This implies that 
	$$\eta(F, a) = \frac{1}{d} \sum_{n=0}^\infty \frac{\sigma(F, f^n(a))}{d^n} = 
		\frac{1}{d^p} \sum_{k=0}^\infty \frac{\sigma(F^p, f^{kp}(a))}{d^{kp}} = \eta(F^p, a).$$
In particular, the irrationality of $\eta(F^p, a)$ (and therefore of the local canonical height for $f^p$, by Lemma \ref{local height}) implies the irrationality of $\lhat_{f,v}(a)$.   

This concludes the proof of Theorem \ref{quadratic}.

\bigskip
\section{Orbit intersections are unlikely}
\label{section diophantine}
We work under the hypotheses of Theorem~\ref{main result application}.

We assume $\OO_f(a)\cap\OO_g(b)$ is infinite and derive a contradiction. In particular, our assumption yields that $a$ (respectively, $b$) is not a preperiodic point for $f$ (respectively, $g$).

Under hypothesis~(2) in Theorem~\ref{main result application}, it follows from 
\cite[Theorem 1.2, Proposition 2.3]{D:stableheight} that 
$$H_1:=\hhat_f(a)>0$$
because $a$ is not preperiodic.  
Set $H_2:=\hhat_g(b)$.  Since the difference between the usual Weil height and the canonical height $\hhat_f$ or $\hhat_g$ is uniformly bounded on $\P^1(\kbar)$, there exists a constant $C>0$ so that 
\begin{equation}
\label{1st equation}
\left|\hhat_f\left(f^m(a)\right) - \hhat_g\left(g^n(b)\right)\right| \leq C,
\end{equation}
for all $m,n\in\N$ such that $f^m(a)=g^n(b)$. Then \eqref{1st equation} yields that 
$$\left|d^m\cdot H_1 - e^n\cdot H_2\right| \leq C,$$
for the infinite collection of pairs of positive integers $(m,n)$ for which $f^m(a)=g^n(b)$.  In particular, the positivity of $H_1$ implies that we also have $H_2 > 0$.  
Furthermore, since $H_1$ and $H_2$ are rational numbers (by hypothesis~(3)), we conclude that there exists a finite set $S$ of rational numbers such that for each $m,n\in\N$ satisfying $f^m(a)=g^n(b)$, we have that
\begin{equation}
\label{2nd equation}
d^m\cdot H_1 - e^n\cdot H_2 \in S.
\end{equation}
By the pigeonhole principle, there exists a rational number $q\in S$ for which there exist infinitely many $(m,n)\in\N\times \N$ with
\begin{equation}
\label{3rd equation}
d^m\cdot H_1 - e^n\cdot H_2 = q.
\end{equation}
This means that the affine line $L$ given by the equation 
\begin{equation}
\label{equation of a line}
H_1\cdot x - H_2\cdot y = q
\end{equation}
contains infinitely many points in common with the finitely generated subgroup $\Gamma$ of $K^*\times K^*$ which consists of all points $\left(d^k,e^\ell\right)$, for $k,\ell\in\Z$. Then Lang's theorem on subgroups of $\bG_m^2$ \cite{Lang}  
yields that  $L \cap (K^*\times K^*)$ is a coset of an algebraic torus, so that $q=0$. Note that when $q=0$, then $L$ is defined by the equation $y=(H_1/H_2)\cdot x$, which is the multiplicative translate by the point $(1, H_1/H_2)$ of the diagonal line $y=x$, which is isomorphic to $\bG_m$, when we see $L$ as a subset of $\bG_m\times \bG_m$. A nonzero value of $q$ does not make $L$ a coset of an algebraic torus.

So with $q=0$, the equation \eqref{3rd equation} reads
\begin{equation}
\label{4th equation}
d^m\cdot H_1 =e^n\cdot H_2
\end{equation}
for infinitely many pairs $(m,n)\in\N\times\N$.  Taking two distinct pairs $(m_1,n_1)$ and $(m_2,n_2)$ satisfying \eqref{4th equation}, and using the fact that $H_1$ and $H_2$ are nonzero, we obtain that 
\begin{equation}
\label{multiplicative dependence d and e}
d^{m_2-m_1}=e^{n_2-n_1}, 
\end{equation}
which contradicts hypothesis~(1) and thus proves that the intersection $\OO_f(a)\cap \OO_g(b)$ must be finite if hypotheses~(1)-(3) are met.

Now, assume that hypothesis~(1) is replaced by the weaker hypothesis~(1'). Then equation~\eqref{multiplicative dependence d and e} yields that $d$ and $e$ must be multiplicative dependent, i.e., there exist coprime positive integers $k$ and $\ell$ such that $d^k=e^\ell$. Furthermore, we get that \eqref{multiplicative dependence d and e} yields the existence of an integer $s$ such that 
\begin{equation}
\label{s 1}
m_2-m_1=ks\text{ and }n_2-n_1=\ell s.
\end{equation}
So, if $\OO_f(a)\cap\OO_g(b)$ is infinite, then \eqref{s 1} yields that there exist infinitely many $s\in\mathbb{N}$ such that $f^{m_1+sk}(a)=g^{n_1+s\ell}(b)$, i.e., the diagonal $\Delta$ of $\mathbb{P}^1\times \mathbb{P}^1$ contains infinitely many points from the orbit of the point $\left(f^{m_1}(a), g^{n_1}(b)\right)$ under the action of the endomorphism $\Phi$ of $\mathbb{P}^1\times\mathbb{P}^1$ given by:
$$(x,y)\mapsto \left(f^k(x), g^\ell(y)\right).$$
Then the Dynamical Mordell-Lang Conjecture (see \cite{GT-JNT}) yields that $\Delta$ is periodic under the action of $\Phi$. This means that there exists $r\in\mathbb{N}$ such that $f^{kr}=g^{\ell r}$, contradicting thus hypothesis~(1'), which completes the proof of Theorem~\ref{main result application}.

\bigskip \bigskip
\def\cprime{$'$}


\bigskip\bigskip

\end{document}